\documentclass[11pt,oneside]{amsart}
\usepackage{amsmath}
\usepackage{amsfonts}
\usepackage{amssymb}
\usepackage{amsthm}
\usepackage{extpfeil}
\usepackage[colorlinks]{hyperref} 
\usepackage{wasysym}
\usepackage{graphicx}
\usepackage{dashrule}
\usepackage[all]{xy}
\usepackage{multirow}
\usepackage{enumitem}
\usepackage{tikz}
\usetikzlibrary{arrows,decorations.pathmorphing,backgrounds,positioning,fit,calc}
\usepackage{tikz-cd}
\usepackage{framed}

\normalfont\upshape

\numberwithin{equation}{section}

\DeclareMathOperator{\GL}{GL}
\DeclareMathOperator{\SL}{SL}

\DeclareMathOperator{\spec}{Spec}

\DeclareMathOperator{\conv}{Conv}
\DeclareMathOperator{\stab}{Stab}

\DeclareMathOperator{\Aut}{Aut}
\DeclareMathOperator{\Lie}{Lie}

\DeclareMathOperator{\quot}{Quot}

\newcommand{\nocontentsline}[3]{}
\newcommand{\tocless}[2]{\bgroup\let\addcontentsline=\nocontentsline#1{#2}\egroup}

\newcommand{\CC}{\mathbb{C}}
\newcommand{\PP}{\mathbb{P}}

\newcommand{\GG}{\mathbb{G}}

\newcommand{\dblslash}{/\! \!/}

\newcommand{\env}{\!
\mathbin{\text{\rotatebox[origin=c]{70}{\scalebox{1.2}{$\approx$}}}} \!}

\newcommand{\inenv}{\dblslash \!_{\circ}}

\newcommand{\ten}{\otimes}

\newcommand{\kk}{\Bbbk}

\newcommand{\ra}{\rightarrow}

\newcommand{\Proj}{{\rm Proj}}
\newcommand{\Spec}{{\rm Spec}}
\newcommand{\nc}{\newcommand}
\nc{\bla}{\phantom{bbbbb}}

\newcommand{\beq}{\begin{equation}}
\newcommand{\eeq}{\end{equation}}
\newcommand{\barr}{\begin{array}}
\newcommand{\earr}{\end{array}}
\newcommand{\beqar}{\begin{eqnarray}}
\newcommand{\eeqar}{\end{eqnarray}}

\newtheorem{thm}{Theorem}[section]
\newtheorem{corollary}[thm]{Corollary}

\newtheorem{prop}[thm]{Proposition}

\newtheorem{exit}[thm]{Example}

\theoremstyle{definition}
\newtheorem{definition}[thm]{Definition}
\newtheorem{rmk}[thm]{Remark}

\newcommand{\RR}{{\mathbb R }}
\nc{\FF}{ {\mathbb F} }
\nc{\HH}{ {\mathbb H} }
\newcommand{\ZZ}{{\mathbb Z }}

\newcommand{\cO}{\mathcal{O}}
\newcommand{\cM}{\mathcal{M}}
\newcommand{\cL}{\mathcal{L}}
\newcommand{\cE}{\mathcal{E}}
\newcommand{\cF}{\mathcal{F}}

\newcommand{\calo}{\mathcal{O}}

\newcommand{\calr}{\mathcal{R}}

\nc{\umax}{{U_{\max}}}

\newcommand{\reg}{\mathrm{reg}}

\newcommand{\hU}{\widehat{U}}

\def\a{\alpha}
\def\b{\beta}

\title{Stratifying quotient stacks and moduli stacks}

\author{Gergely B\'erczi, Victoria Hoskins,  Frances Kirwan}

\address[B\'erczi]{Departement Mathematik, ETH, 
 Switzerland}
\email{gergely.berczi@math.ethz.ch}

\address[Hoskins]{Fachbereich Mathematik und Informatik\\
Freie Universit\"{a}t Berlin\\Germany} 
\email{hoskins@zedat.fu-berlin.de}
\address[Kirwan]{Mathematical Institute, Oxford University, UK} 
\email{kirwan@maths.ox.ac.uk}

\thanks{V.H. is supported by the Excellence Initiative of the DFG at the Freie Universit\"{a}t Berlin}

\begin{document}

\begin{abstract}
Recent results in geometric invariant theory (GIT) for non-reductive linear algebraic group actions allow us to stratify quotient stacks of the form $[X/H]$, where $X$ is a projective scheme and $H$ is a linear algebraic group with internally graded unipotent radical acting linearly on $X$, in such a way that each stratum $[S/H]$ has a geometric quotient $S/H$. This leads to stratifications of moduli stacks (for example, sheaves over a projective scheme) such that each stratum has a coarse moduli space.
\end{abstract}

\maketitle

\section{Introduction}\label{sec intro}

Let $H = U \rtimes R$ be a linear algebraic group over an algebraically closed field of characteristic 0 with internally graded unipotent radical $U$; that is, the Levi subgroup $R$ of $H$ has a central one-parameter subgroup (1-PS) $\lambda: \GG_m \to R$ which acts on $\Lie U$ with all weights strictly positive. Of course any reductive group $G$ has this form with both $U$ and the central one-parameter subgroup $\lambda$ being trivial. Suppose that $H$ acts linearly on an irreducible projective scheme $X$ with respect to an ample line bundle $L$ over $X$. The aim of this paper is to describe a stratification of the quotient stack $[X/H]$ such that each stratum $[S/H]$ (where $S$ is an $H$-invariant quasi-projective subscheme of $X$) has a geometric quotient $S/H$. When $H = R$ is reductive this stratification refines the Hesselink--Kempf--Kirwan--Ness (HKKN) stratification associated to the linear action on $X$ (\textit{cf.}\ \cite{Hess,Kempf,K,Ness}). Potential applications of this construction include moduli stacks which can be filtered by quotient stacks with compatible linearisations; for example, it can be applied to moduli of  sheaves of fixed Harder--Narasimhan type over a projective scheme \cite{behjk16}, and moduli of unstable projective curves \cite{jj}. 

When $H$ is reductive Mumford's geometric invariant theory (GIT) \cite{GIT} allows us to find open subschemes $X^s \subseteq X^{ss}$ of $X$, the stable and semistable loci, such that $X^s$ has a geometric quotient $X^s/H$ and $X^{ss}$ has a good quotient 
$$\phi: X^{ss} \to X/\!/H =  \Proj \left(\bigoplus_{m \geqslant 0} H^0(X,L^{\otimes m})^H \right).$$
Here $X^s= X^{s,H} =X^{s,H,\mathcal{L}}$ and $X^{ss}= X^{ss,H} =X^{ss,H,\mathcal{L}}$  depend on the choice of linearisation $\mathcal{L}$ (that is, the ample line bundle $L$ and the lift of the action of $H$ to an action on $L$) and the GIT quotient $X/\!/H = X/\!/\!_{\mathcal{L}}H$ is a projective scheme with $X^s/H$ as an open subscheme. Moreover when $x,y \in X^{ss}$ then $\phi(x) = \phi(y)$ if and only if the closures of the $H$-orbits of $x$ and $y$ meet in $X^{ss}$. The Hilbert--Mumford criteria allow us to determine the stable and semistable loci in a simple way without needing to understand the algebra of invariants $\bigoplus_{m \geqslant 0} H^0(X,L^{\otimes m})^H$.

The best situation occurs when $X^{ss} = X^s \neq \emptyset$; then $X/\!/H = X^s/H$ is both a projective scheme and a geometric quotient of the open subscheme $X^s$ of $X$. More generally if $X^s \neq \emptyset$ then the projective completion $X/\!/H$ of the geometric quotient $X^s/H$ has a \lq partial desingularisation' $\widetilde{X}/\!/H = \widetilde{X}^{ss}/H$ where $\psi: \widetilde{X}^{ss} \to X^{ss}$ is obtained as follows \cite{K2}.
If $X^{ss} \neq X^s$ then there exists $x \in X^{ss}$ whose stabiliser in $H$ is reductive of dimension strictly bigger than 0. To construct $\widetilde{X}^{ss}$ we first blow up $X^{ss}$  along its closed subscheme where the stabiliser has maximal dimension in $X^{ss}$ (such stabilisers are always reductive) or equivalently blow up  $X$ along the closure of this subscheme in $X$. We  then remove the complement of the semistable locus for a small ample perturbation of the pullback linearisation. The maximal dimension of a stabiliser in this new semistable locus is strictly less than in $X^{ss}$. When $X^s \neq \emptyset$, repeating this process finitely many times leads to $\widetilde{X}^{ss} = \widetilde{X}^s \neq \emptyset$ and the partial desingularisation $\widetilde{X}/\!/H = \widetilde{X}^{s}/H$.

When $H$ is non-reductive, then the graded algebra $\bigoplus_{m \geqslant 0} H^0(X,L^{\otimes m})^H$ is not necessarily finitely generated and in general the attractive properties of Mumford's GIT fail \cite{BDHK}. However when the unipotent radical $U$ of $H = U \rtimes R$ is graded in the sense described above 
by a central 1-PS $\lambda:\GG_m \to R$ of the Levi subgroup $R$, then after twisting the linearisation by an appropriate rational character, so that it becomes \lq graded' itself in the sense of \cite{BDK}, some of the desirable properties of classical GIT still hold \cite{BDHK2,BDHK3}. More precisely, we first quotient by the linear action of the graded unipotent group $\hU := U \rtimes \lambda(\GG_m)$ using the results of \cite{BDHK2,BDHK3} described in $\S$\ref{GIT_int_gr}, then we quotient by the residual action of the reductive group $H/\hU \cong R/\lambda(\GG_m)$. In the best case, when the $U$-action is free on a certain open subscheme $X^0_{\min}$ of $X$ (\textit{cf.}\ the condition \eqref{star} in Definition \ref{cond star} and Theorem \ref{mainthm}), one can construct a geometric quotient of an open subscheme of \lq stable points' for the $\hU$-action such that the quotient is projective and this stable set has a Hilbert--Mumford type description. If the $U$-action has positive dimensional stabilisers generically, one can conclude the same results if we assume a weaker condition \eqref{star3} (\textit{cf.}\ Theorem \ref{Uh3}). Even when this weaker condition fails, one can perform an iterated sequence of blow-ups of $X$ along $H$-invariant subschemes to obtain $\psi  : \widetilde{X} \ra X$ such that $\widetilde{X}$ has an induced linear $H$-action satisfying \eqref{star3}. Hence, there is a projective and geometric $\hU$-quotient of an open subscheme of $\widetilde{X}$ that contains as an open subscheme a geometric $\hU$-quotient of an open subscheme of $X$, as $\psi$ is an isomorphism away from the exceptional divisor.

Now suppose that  $G$ is a reductive group  acting linearly on a projective scheme $X$ with respect to an ample line bundle $L$. Associated to this linear $G$-action and an invariant inner product on $\Lie G$, there is a stratification (the \lq HKKN stratification' \textit{cf.}\ \cite{Hess,Kempf,K,Ness})
$$ X = \bigsqcup_{\beta \in \mathcal{B}} S_\beta$$ of $X$ by locally closed subschemes $S_\beta$, 
indexed by a partially ordered finite set $\mathcal{B}$, such that 
\begin{enumerate}
\item if $X^{ss} \neq \emptyset$, then $\mathcal{B}$ has a minimal element $0$ such that $S_0 = X^{ss}$,
\item for $\beta  \in \mathcal{B}$, the closure of $S_\beta$ is contained in $\bigcup_{\beta' \geqslant \beta} S_{\beta'}$, and
\item for $\beta \in \mathcal{B}$, we have  $S_\beta \cong G \times^{P_\beta} Y_\beta^{ss}$, where $G \times^{P_\beta} Y_\beta^{ss}$ is the quotient of $G \times Y_\beta^{ss}$ by the diagonal action of a parabolic subgroup $P_\beta$ of $G$ acting on the right on $G$ and on the left on a $P_\beta$-invariant locally closed subscheme $Y_\beta^{ss}$ of $X$.
\end{enumerate}
\noindent In fact, $Y_\beta^{ss}$ is an open subscheme of a projective subscheme $\overline{Y}_\beta$ of $X$ that is determined by the action of a Levi subgroup $L_\beta$ of $P_\beta$ with respect to  the restriction of the $G$-linearisation $L \to X$ to the $P_\beta$-action on $\overline{Y}_\beta$ twisted by a rational character $\chi_\beta$ of $P_\beta$. The index $\beta$ determines a central (rational) 1-PS $\lambda_\beta : \GG_m \ra L_\beta$ and $\chi_\beta$ is the corresponding rational character, where the choice of invariant inner product allows us to identify characters and co-characters of a fixed maximal torus (\textit{cf.}\ Remark \ref{inner prod}). Furthermore $P_\beta$ is the parabolic subgroup  $P(\lambda_\beta)$ determined by the 1-PS $\lambda_\beta$, which grades the unipotent radical $U_\beta$ of $P_\beta$.  Thus by Property (3) above, to construct a $G$-quotient of (an open subset of) an unstable stratum $S_\beta$, we can study the linear $P_\beta$-action on $\overline{Y}_\beta$  and apply the  results described above  for the action of $\hU := U_\beta \rtimes \lambda_\beta(\GG_m)$.

The $G$-action on the stratum $S_\beta$ has a categorical quotient $Z_\beta/\!/L_\beta$ induced by the morphism 
\[ p_\beta: Y_\beta^{ss} \to Z_\beta^{ss} = Z_\beta^{ss,L_\beta/\lambda_\beta(\GG_m)} \quad \quad  y \mapsto p_\beta(y):=\lim_{t\to 0} \lambda_\beta(t)y,\]
where $Z_\beta$ is the union of those connected components of the fixed point set for the action of $\lambda_\beta(\GG_m)$ on $X$ over which $\lambda_\beta$ acts on the fibres of $L$ with weight given by the restriction of $\beta$. However this categorical quotient is in general far from being a geometric quotient; 
 it identifies $y$ with $p_\beta(y)$, which lies in the orbit closure but typically not the orbit of $y$.

In this article we will show that, applying the blow-up sequence needed to construct a quotient by an action of a linear algebraic group with internally graded unipotent radical to the $P_\beta$-action on the projective subscheme $\overline{Y}_\beta$ of $X$, we can refine the stratification $\{S_\beta | \beta \in \mathcal{B} \}$ to obtain a stratification of $X$ such that each stratum is a $G$-invariant quasi-projective subscheme of $X$ with a geometric quotient by the action of $G$. This refined stratification is a further refinement of the construction described in \cite{refinement}. The quotient stack $[X/G]$ has an induced stratification $\{\Sigma_\gamma | \gamma \in \Gamma \}$ such that each stratum $\Sigma_\gamma$ has the form 
$$ \Sigma_\gamma \cong [W_\gamma/H_\gamma]$$
where $W_\gamma$ is a quasi-projective subscheme of $X$ acted on by a linear algebraic subgroup $H_\gamma$ of $G$ with internally graded unipotent radical, and this action has a geometric quotient $W_\gamma/H_\gamma$. Moreover under appropriate hypotheses (involving condition \eqref{star3}) for the actions of the subgroups $H_\gamma$, the geometric quotients $W_\gamma/H_\gamma$ are themselves projective.

This will follow from the following theorem, which is proved in $\S$\ref{sec strat quot stacks}.

\begin{thm} \label{mainthmstrat} Let $H = U \rtimes R$ be a linear algebraic group with internally graded unipotent radical $U$ acting on a projective scheme $X$ over an algebraically closed field $\kk$ of characteristic 0 with respect to an ample linearisation and fix an invariant inner product on $\Lie R$. Then $X$ has a stratification $\{ \mathcal{S}_\gamma | \gamma \in \Gamma \}$ induced by the linearisation $\mathcal{L}$ and grading $\lambda:\GG_m \to R$ for the action of $H$ on $X$, such that the following properties hold.
\begin{enumerate}[label=\emph{\roman*})]
\item The index set $\Gamma$ is finite and partially ordered such that for all $\gamma \in \Gamma$, we have
$$\overline{\mathcal{S}_\gamma} \subseteq \mathcal{S}_\gamma \cup \bigcup_{\delta \in \Gamma, \delta > \gamma} \mathcal{S}_\delta.$$
\item  Each $\mathcal{S}_\gamma$ is a $H$-invariant quasi-projective subscheme of $X$ with a geometric quotient $\mathcal{S}_\gamma /H$.
\item  If $Y$ is an $H$-invariant projective subscheme of $X$ then the stratification $\{ \mathcal{S}^Y_\gamma | \gamma \in \Gamma^Y \}$ of $Y$ induced by the restriction $\mathcal{L}|_Y$ of the linearisation $\mathcal{L}$ to Y is (up to taking connected components) the restriction to $Y$ of the stratification $\{ \mathcal{S}_\gamma | \gamma \in \Gamma \}$ of $X$, so that there is a map of indexing sets $\phi_Y:\Gamma^Y \to \Gamma$ such that if $\gamma \in \Gamma^Y$ then $\mathcal{S}^Y_\gamma$ is a connected component of $\mathcal{S}_{\phi_Y(\gamma)} \cap Y$;
\end{enumerate}
\noindent Moreover, if $H=G$ is reductive, then this stratification satisfies the following additional properties.
\begin{enumerate}[label=\emph{\roman*}),resume]
\item The stratification $\{ \mathcal{S}_\gamma | \gamma \in \Gamma \}$ is a refinement of the HKKN stratification $\{S_\beta | \beta \in \mathcal{B} \}$ for the linearisation $\mathcal{L}$ (\textit{cf.}\ $\S$\ref{stratification}).
\item If $\beta \in \mathcal{B}$ (which we recall determines a 1-PS $\lambda_\beta: \GG_m \to P_\beta \leq G$) satisfies
\begin{equation}\label{dagger}
 x \in Z_{\beta}^{ss} \quad \Rightarrow  \quad \dim (\stab_G(x) / \lambda_\beta(\GG_m)) = 0, \tag{$\dagger$}
\end{equation}
then the connected components of $GZ_\beta^{ss}$ and $S_\beta \setminus GZ_\beta^{ss}$ (if these are nonempty) are strata in the refined stratification $\{ \mathcal{S}_\gamma | \gamma \in \Gamma \}$.
\end{enumerate}
\end{thm}

As a consequence, we obtain a stratification of the quotient stack $[X/H]$ by locally closed substacks, each of which admits a coarse moduli space (\textit{cf.}\ Corollary \ref{cor cms}). Inspired by the reductive GIT notion of a good quotient, Alper introduces a notion of a good moduli space for a stack \cite{alper}. However, in general the strata appearing in this stratification of $[X/H]$ will not admit good moduli spaces, because a necessary condition for a stack to admit a good moduli space is that its closed points have reductive stabiliser groups (\textit{cf.}\ \cite[Proposition 12.14]{alper}). In general  (even when $H = G$ is reductive) the points in the strata of $[X/H]$ will have non-reductive stabiliser groups. 

If $H =G$ is reductive, then a stacky version of the HKKN stratification has been studied by Halpern-Leister, and by abstracting this concept he obtains a notion of a $\Theta$-stratification \cite{HL}. Indeed, the linearisation of the $G$-action on $X$ and the choice of invariant inner product is precisely the data required to construct a $\Theta$-stratification of $[X/G]$, and this $\Theta$-stratification is the stratification $\{ [\mathcal{S}_\beta/G] : \beta \in \mathcal{B} \}$ obtained from the HKKN stratification of $X$. The  stratification described above thus refines this $\Theta$-stratification without depending on any additional data.

Since the construction of the refined stratification involves studying the blow-up procedures used in partial desingularisations of reductive GIT quotients \cite{K2} and for constructing geometric quotients by linear algebraic groups with internally graded unipotent radical \cite{BDHK3}, one can ask how this compares with the stack-theoretic blow-up constructions. The ideas in \cite{K2} have been generalised to stacks by Edidin and Rydh \cite{ER} to show that for a smooth Artin stack $\mathfrak{X}$ admitting a stable good moduli space, there is a sequence of birational morphisms of smooth Artin stacks $\mathfrak{X}_n \ra \cdots \mathfrak{X}_1 \ra \mathfrak{X}$ such that the good moduli space of $\mathfrak{X}_n$ is an algebraic space with only tame quotient singularities and is a partial desingularisation of the good moduli space of $\mathfrak{X}$. However, for $H$ non-reductive, it is often the case that $[X/H]$ will not be a good moduli space, and so one cannot apply this result. 

The picture  provided by Theorem \ref{mainthmstrat} has potential applications to moduli stacks which are filtered by quotient stacks, and to the construction of moduli spaces of \lq unstable' objects (for example, moduli of sheaves over a projective scheme \cite{behjk16} or moduli of projective curves \cite{jj}). Suppose that $\mathcal{M}$ is a moduli stack which can be expressed as an increasing union 
$$\mathcal{M} = \bigcup_{n \geqslant 0} \mathcal{U}_n $$
of open substacks of the form
$$ \mathcal{U}_n \cong [V_n/G_n]$$
where $[V_n/G_n]$ is the quotient stack associated to a linear action on a quasi-projective scheme $V_n$ by a group $G_n$ which is reductive (or more generally has internally graded unipotent radical). We can look for suitable \lq stability conditions' on $\mathcal{M}$: linearisations $(\mathcal{L}_n)_{n\geqslant 0}$ for the actions of $G_n$ on projective completions $\overline{V_n}$ of $V_n$ and invariant inner products on $\Lie G_n$ which are compatible in the sense that the stratification induced by $\mathcal{L}_n$ on $[V_n/G_n]$ restricts to the stratification induced by $\mathcal{L}_{m}$ on $[V_{m}/G_{m}]$ when $n>m$. 

This situation arises for sheaves over a projective scheme \cite{hok12,behjk16}, for example, and also for projective curves \cite{jj}, and we obtain a stratification  $\{\Sigma_\gamma | \gamma \in \Gamma \}$ of the stack $\mathcal{M}$ such that each stratum $\Sigma_\gamma$ is isomorphic to a quotient stack 
$ [W_\gamma/H_\gamma]$,
where $W_\gamma$ is  quasi-projective acted on by a linear algebraic group $H_\gamma$  with internally graded unipotent radical, and there is a geometric quotient $W_\gamma/H_\gamma$ which is a coarse moduli space for $\Sigma_\gamma$. The geometric quotient $W_\gamma/H_\gamma$ will be projective if semistability coincides with stability in an appropriate sense for the action of $H_\gamma$ on a suitable projective completion of $W_\gamma$ with respect to an induced linearisation.

The layout of this article is as follows. In $\S$\ref{sec GIT}, we will review classical and non-reductive GIT, describing how to construct quotients by actions of linear algebraic groups with internally graded unipotent radical. The heart of the paper is $\S$\ref{sec strat quot stacks}, in which we describe how to stratify a quotient stack $[X/H]$ into strata $\Sigma_\gamma = [W_\gamma/H_\gamma]$ where the action of $H_\gamma$ on $W_\gamma$ has a geometric quotient $W_\gamma/H_\gamma$. The argument is an inductive one, so the assumption on $X$ and $H$ is that $H$ is a linear algebraic group with internally graded unipotent radical, and $X$ is a projective scheme which has an amply linearised action of $H$. In $\S$\ref{sec filt stack}, this construction is applied to stacks which are suitably filtered by quotient stacks, and $\S$\ref{sec appl} contains a brief discussion of examples including moduli of unstable curves and moduli of sheaves of given Harder--Narasimhan type over a fixed projective scheme. 

We would like to thank Brent Doran, Daniel Halpern-Leistner, Eloise Hamilton and Joshua Jackson for helpful discussions about this material.

\section{Classical and non-reductive GIT}\label{sec GIT}

\subsection{Classical GIT for reductive groups}\label{sec GIT classical}

In Mumford's GIT \cite{GIT}, a linearisation of an action of a reductive group $G$ on a projective scheme $X$ over an algebraically closed field $\kk$ of characteristic $0$ is given by a line bundle $L$ (which we will always assume to be ample) on $X$ and a lift of the action to $L$. Since $G$ is reductive, the algebra of $G$-invariant sections ${\widehat{\cO}}_L(X)^G$ is finitely generated as a graded algebra with associated projective scheme  
$X/\!/G = \mbox{Proj}({\widehat{\calo}}_L(X)^G)$. 
$$\begin{array}{ccccl}
(X,L) & \leadsto 
 & {\widehat{\cO}}_L(X)&:=&  \bigoplus_{k= 0}^{\infty} H^0(X,L^{\otimes k})\\
| &&&& \\
| & & \bigcup \!| & &\\
\downarrow & & & & \\
X/\!/G & {\reflectbox{\ensuremath{\leadsto}}}   
 & {\widehat{\cO}}_L(X)^G & & \mbox{ algebra of invariants. }
\end{array}$$
The inclusion of ${\widehat{\cO}}_L(X)^G$ in ${\widehat{\cO}}_L(X)$ determines a rational map $X - - \rightarrow X/\!/G$ which fits into a diagram
$$\begin{array}{rcccl}
 & X & - - \rightarrow & X/\!/G & \mbox{ projective}\\
 & \bigcup  & & || & \\
\mbox{semistable} & X^{ss} & \xtwoheadrightarrow{\quad} & X/\!/G & \\
 & \bigcup  & & \bigcup & \mbox{open}\\
\mbox{stable} & X^s & \xtwoheadrightarrow{\quad} & X^s/G &
\end{array} $$     
where $X^s$ (the stable locus) and $X^{ss}$ (the semistable locus) are open subschemes of $X$, there is a geometric quotient $X^s/G$ for the action of $G$ on $X^s$, the GIT quotient ${X}/\!/G$ is a good quotient for the action of $G$ on $X^{ss}$ via the $G$-invariant surjective morphism $\phi: X^{ss} \to X/\!/G$, and  
$$\phi(x) = \phi( y) \Leftrightarrow \overline{Gx} \cap \overline{Gy} \cap X^{ss} \neq \emptyset.$$
The semistable and stable loci $X^{ss}$ and $X^s$ of $X$  are characterised by the following properties (\cite[Chapter 2]{GIT}, \cite{New}). 

\begin{prop} [Hilbert--Mumford criteria for reductive groups]
Let $T$ be a maximal torus of $G$.
\begin{enumerate}
\item A point $x \in X$ is semistable (respectively stable) for the $G$-action on $X$ if and only if for every
$g\in G$ the point $gx$ is semistable (respectively stable) for the $T$-action.
\item  A point $x \in X \subset \PP^n$ with homogeneous coordinates $[x_0:\ldots:x_n]$ is semistable (respectively stable) for a diagonal $T$-action on $\PP^n$ with weights $\a_0, \ldots, \a_n$ if and only if 
\[ 0 \in \conv \{\a_i :x_i \neq 0\} \]
(respectively $0$  is contained in the interior of this convex hull).
\end{enumerate}
\end{prop}

\subsubsection{The HKKN stratification} \label{stratification}
Associated to the linear action of $G$ on $X$ and an invariant inner product on the Lie algebra of $G$, there is a stratification (the \lq HKKN stratification', which  in the case $\kk = \CC$ is the Morse stratification for the norm-square of an associated moment map \cite{Hess, Kempf, K, Ness}).

\begin{rmk}\label{inner prod}
Let us clarify what is meant by this invariant inner product, whose associated norm we denote by $|| - ||$. If $\kk = \CC$, then $G$ is the complexification of its maximal compact group $K$; then the Lie algebra of $K$ is a real vector space, and we choose an inner product on this Lie algebra that is invariant under the adjoint action of $K$. In fact, we will also assume that we fix a maximal compact torus $T_c \subset K$ such that the inner product is integral on the co-character lattice $X_*(T_c) \subset \Lie K$. For an arbitrary algebraically closed field $\kk$ of characteristic zero, one can fix a maximal torus $T$ of $G$ and choose an inner product on the co-character space $X_*(T) \otimes_\ZZ \RR$ that is invariant for the Weyl group of $T$ and is integral on the co-character lattice (for example, see \cite[$\S$2]{H}). Then this inner product gives an identification between characters and co-characters (i.e. 1-PSs) of $T$.
\end{rmk}

The HKKN stratification associated to the action of $G$ on $X$ with respect to $L$ and the norm $||-||$ is a stratification
$$ X = \bigsqcup_{\beta \in \mathcal{B}} S_\beta$$ of $X$ by locally closed subschemes $S_\beta$, 
indexed by a partially ordered finite subset $\mathcal{B}$ of a rational elements in a positive Weyl chamber for the reductive group $G$, with the following properties.
\begin{enumerate}
\item If $0 \in \mathcal{B}$, then this is the minimal element and $S_0 = X^{ss}$.
\end{enumerate}
Moreover, for each $\beta \in \mathcal{B}$, we additionally have the following properties.
\begin{enumerate}[resume]
\item the closure of $S_\beta$ is contained in $\bigcup_{\beta' \geqslant \beta} S_{\beta'}$  where $\gamma \geqslant \beta \mbox{ if and only if } \gamma = \beta \mbox{ or } |\!|\gamma|\!| > |\!|\beta|\!|$;
\item \label{strata str} $S_\beta \cong G \times^{P_\beta} Y_\beta^{ss}:= (G \times Y_\beta^{ss})/P_\beta$ where this quotient is of the diagonal action of $P_\beta$ on the right on $G$ and on the left on $Y^{ss}_\beta$.
\end{enumerate}
\noindent Here $P_\beta$ is a parabolic subgroup of $G$ which acts on a locally closed subscheme $Y^{ss}_\beta$ of $X$. 

More precisely, $\beta \in \mathcal{B}$ determines a (rational) 1-PS $\lambda_\beta : \GG_m \ra G$ and an associated parabolic subgroup $P_\beta = P(\lambda_\beta) = U_\beta \rtimes L_\beta$ with Levi subgroup $L_\beta = \stab_G(\beta)$ such that the conjugation action of $\lambda_\beta(\GG_m)$ on $\Lie U_\beta$ has strictly positive weights; thus $P_\beta$ has internally graded unipotent radical. Let $Z_\beta$ be the union of components in the fixed locus $X^{\lambda_\beta(\GG_m)}$ on which this 1-PS acts on the fibres of $L$ with weight given by the restriction of $\beta$, and let $Y_\beta \subset X$ be the subscheme of points $x \in X$ such that $\lim_{t \ra 0} \lambda_\beta(t) y \in Z_\beta$; thus there is a retraction 
\[p_\beta : Y_\beta \ra Z_\beta \quad \quad y \mapsto p_\beta(y) := \lim_{t \ra 0} \lambda_\beta(t) y\] 
which is equivariant with respect to the quotient homomorphism $q_\beta : P_\beta \ra L_\beta$ obtained by identifying $L_\beta$ with $P_\beta/U_\beta$.  Let $\cL_\beta$ denote the restriction of the $G$-linearisation $\cL$ on $X$ to the $P_\beta$-action on $\overline{Y}_\beta$ twisted by the (rational) character $\chi_\beta$ of $P_\beta$ corresponding to the 1-PS $\lambda_\beta$ (via the norm $|| - ||$). We also let $\cL_\beta$ denote the restriction of this linearisation to the $L_\beta$-action on $Z_\beta$. Then $Y_\beta^{ss}$ (respectively $Z_\beta^{ss}$) is the semistable locus for the $P_\beta$-action on $\overline{Y}_\beta$ (respectively the $L_\beta/\lambda_\beta(\GG_m)$-action on $Z_\beta$) linearised by the twisted linearisation $\cL_\beta$; furthermore, $Y_\beta^{ss} = p_\beta^{-1}(Z_\beta^{ss})$.

Finally, we make the following observation about quotienting the unstable strata (\textit{cf.}\ \cite{hok12}).

\begin{rmk}
The $G$-action on $S_\beta$ has a categorical quotient $Z_\beta/\!/L_\beta$ induced by the map $p_\beta: Y_\beta^{ss} \to Z_\beta^{ss}$. In general, this quotient is far from being a geometric quotient (even after restriction to the pre-image of any nonempty open subscheme of   $Z_\beta/\!/L_\beta$), as $y \in Y_\beta^{ss}$ is identified with $p_\beta(y) \in \overline{G y}$.
\end{rmk}

By \eqref{strata str}, constructing a quotient of the $G$-action on a $G$-invariant open subset of $S_\beta$ is equivalent to constructing a $P_\beta$-quotient of a $P_\beta$-invariant open subset of $Y_\beta^{ss}$ (or its closure); the latter perspective will lead to a geometric quotient 
 by using GIT for the non-reductive group $P_\beta$, whose unipotent radical $U_\beta$ is internally graded by $\lambda_\beta$ (\textit{cf.}\ $\S$\ref{GIT_int_gr}).

\subsubsection{Partial desingularisations of reductive GIT quotients}\label{sec part desing}

The geometric quotient $X^s/G$ has at most orbifold singularities when $X$ is nonsingular,  since the stabiliser subgroups of stable points are finite subgroups of $G$. If $X^{ss} \neq X^s \neq \emptyset$, the singularities of $X/\!/G$ are typically more severe even when $X$ is itself nonsingular, but $X/\!/G$ has a \lq partial desingularisation'  $\widetilde{X}/\!/G$ which is also a projective completion of $X^s/G$ and  is itself a geometric quotient 
$$\widetilde{X}/\!/G = \widetilde{X}^{ss}/G$$
by $G$ of an open subscheme $\widetilde{X}^{ss} = \widetilde{X}^s$ of a $G$-equivariant blow-up $\widetilde{X}$ of $X$ \cite{K2}.  Here $\widetilde{X}^{ss}$ is obtained from ${X}^{ss}$ by successively blowing up along the subschemes of semistable points stabilised by reductive subgroups of $G$ of maximal dimension and removing the complement of the semistable locus from the resulting blow-up. 

For the construction of the partial desingularisation $\widetilde{X}/\!/G$ in \cite{K2}, it is assumed that $X^s \neq \emptyset$. There exist semistable points of $X$ which are not stable if and only if there exists a non-trivial connected reductive subgroup of $G$ fixing a semistable point. Let $r>0$ be the maximal dimension of a reductive subgroup of $G$ fixing a point of $X^{ss}$ and let $\calr(r)$ be a set of representatives of conjugacy 
classes of all connected reductive subgroups $R$ of 
dimension $r$ in $G$ such that 
\[ Z^{ss}_{R} := \{ x \in X^{ss} :  Rx = x\} \]
is non-empty. Then
\[
Z^{ss}_{\calr(r)}:=\bigcup_{R \in \calr(r)} GZ^{ss}_{R}
\]
is a disjoint union of closed $G$-invariant subschemes of $X^{ss}$. The action of 
$G$ on $X^{ss}$ lifts to an action on the blow-up $X_{(1)}$ of 
$X^{ss}$ along $Z^{ss}_{\calr(r)}$, and this action can be linearised so that the complement 
of $X_{(1)}^{ss}$ in $X_{(1)}$ is the proper transform of the 
closed subscheme $\pi^{-1}(\pi(GZ_R^{ss}))$ of $X^{ss}$ where $\pi:X^{ss} \to X/\!/G$ is the quotient 
map (see \cite[7.17]{K2}). The $G$-linearisation on $X_{(1)}$ used here is (a tensor power of) the pullback of the ample line bundle $L$ on $X$ along $\psi_{(1)}:X_{(1)} \to X$ perturbed by a sufficiently small multiple of the exceptional divisor $E_{(1)}$; then, if the perturbation is sufficiently small, we have
$$   \psi_{(1)}^{-1}(X^{s}) \subseteq  X_{(1)}^{s} \subseteq  X_{(1)}^{ss} \subseteq \psi_{(1)}^{-1}(X^{ss}) = X_{(1)},$$
and the stable and semistable loci $X_{(1)}^{s}$ and $X_{(1)}^{ss}$ will be independent of the choice of perturbation.
Moreover, no point $x \in X_{(1)}^{ss}$ is fixed by a 
reductive subgroup of $G$ of dimension at least $r$, and  $x \in X_{(1)}^{ss} $ 
is fixed by a reductive subgroup $R$ of 
dimension less than $r$ in $G$ if and only if it belongs to the proper transform of the closed
subscheme $Z_R^{ss}$ of $X^{ss}$.

\begin{rmk} \label{cfK2} In \cite{K2}, $X$ itself is blown up along the closure $\overline{Z^{ss}_{\calr(r)}}$ of $Z^{ss}_{\calr(r)}$ in $X$ (or in a projective completion of $X^{ss}$ with a $G$-equivariant morphism to $X$ which is an isomorphism over $X^{ss}$). This gives a projective scheme $\overline{X}_{(1)}$ and blow-down map $\overline{\psi_{(1)}}: \overline{X}_{(1)} \to X$ restricting to $\psi_{(1)}:X_{(1)} \to X$ where $(\overline{\psi_{(1)}})^{-1}(X^{ss}) = X_{(1)}$.  For a sufficiently small perturbation of the pullback to $\overline{X}_{(1)}$ of the linearisation on $X$,
we have
$   (\overline{\psi_{(1)}})^{-1}(X^{s}) \subseteq \overline{X}_{(1)}^{s} \subseteq \overline{X}_{(1)}^{ss} \subseteq (\overline{\psi_{(1)}})^{-1}(X^{ss}) = X_{(1)},$
and moreover the restriction of the linearisation to $X_{(1)}$ is obtained from the pullback of $L$ by perturbing  by a sufficiently small multiple of the exceptional divisor $E_{(1)}$.
\end{rmk}

If $r>1$, we can apply the same procedure to $X_{(1)}^{ss}$ to obtain $X_{(2)}^{ss}$ such that no point of $X_{(2)}^{ss}$ is fixed by a reductive subgroup of $G$ of dimension at least $r-1$. If $X^s \neq \emptyset$ then repeating this process at most $r$ times gives us $\psi: \widetilde{X}^{ss} \to X^{ss}$ such that $\psi$ is an isomorphism over $X^s$ and no positive-dimensional reductive subgroup of $G$ fixes a point of $\widetilde{X}^{ss}$.  The partial desingularisation $\widetilde{X}/\!/G = \widetilde{X}^{ss}/G$ can be obtained by blowing up $X/\!/G$ along the proper transforms of $\pi(GZ_R^{ss}) \subset X/\!/G$ in decreasing order of the dimension of $R$.

\begin{rmk} \label{rem:specialcases} Suppose for simplicity that $X$ is irreducible.\
\begin{enumerate}[label=\emph{\alph*})]
\item If $X^s \neq \emptyset$, then this is the situation considered in \cite{K2} and the partial desingularisation construction is described above. If $X^{ss} = X^{s}$, then $\widetilde{X} = X$.
\item \label{b} If $X^{ss} = \emptyset$, then there is an unstable stratum $S_\beta$ with $\beta \neq 0$ in the HKKN stratification (\textit{cf.}\ $\S$\ref{stratification}) which is a non-empty open subscheme of $X$, and thus when $X$ is irreducible  $X = \overline{S_\beta}$. Then constructing a quotient of a non-empty open subscheme of $X$ reduces to non-reductive GIT for the action of the parabolic subgroup $P_\beta$ on $\overline{Y_\beta}$ as described in $\S$\ref{GIT_int_gr} below, where a blow-up sequence may also need to be performed.
\item  If $X^s = \emptyset \neq X^{ss}$  then the  partial desingularisation construction can be applied to $X^{ss}$, and there are different ways in which it can terminate.
\begin{enumerate}[label=(\emph{\roman*})]
\item \label{i} If $X^{ss} = GZ_R^{ss} \cong  G \times^{N_R} Z_R^{ss}$ for a positive-dimensional connected reductive subgroup $R$ of $G$ with normaliser $N_R$ in $G$, then $N_R$ and $N_R/R$ are also reductive, and
$$X/\!/G \cong Z_R /\!/ N_R \cong Z_R /\!/ (N_R/R)$$
where $Z_R$ is the closed subscheme of $X$ which is the fixed point set for the action of $R$. Then we can apply induction on the dimension of $G$ to study this case. 

\item If $GZ_R^{ss} \neq X^{ss}$ for each positive-dimensional connected reductive subgroup $R$ of $G$, then we can perform the first blow-up in the partial desingularisation construction to obtain $\psi_{(1)}: X_{(1)} \to X^{ss}$ such that $X_{(1)}^{ss} \subseteq X_{(1)}$ and $X_{(1)}^s = \emptyset$ as above (as $X^s_{(1)}$ is open and $X^s_{(1)} \setminus E_{(1)} \cong X^s = \emptyset$, where $E_{(1)}$ is the exceptional divisor). If $X_{(1)}^{ss} = \emptyset$, then $X_{(1)}$ has a dense open stratum $S_{(1),\beta}$ for $\beta \neq 0$ as in Case \ref{b}. If we have $X_{(1)}^{ss} = GZ_{(1),R}^{ss}$ for a positive-dimensional connected reductive subgroup $R$ of $G$, where $Z_{(1),R}^{ss} = \{ x \in X_{(1)}^{ss} : Rx=x \}$, then we proceed as in Case \ref{i} above. Otherwise we can repeat the process, until it terminates in one of these two ways.
\end{enumerate}
\end{enumerate}
\end{rmk}

\subsection{GIT for non-reductive groups}\label{sec non-red}

Now suppose that $X$ is a projective scheme over an algebraically closed field $\kk$ of characteristic 0 and let $H$ be a linear algebraic group, with unipotent radical $U$, acting on $X$ with respect to an ample linearisation $L$.  

\begin{definition} (Semistability for the unipotent group \textit{cf.}\ \cite[$\S$4]{DK} and \cite[5.3.7]{DK}). \label{defnssetc}
For an invariant section $f \in I = \bigcup_{m>0} H^0(X,L^{\otimes m})^U$, let $X_f$ be the $U$-invariant affine open subset of $X$ on which $f$ does not vanish, and let ${\calo}(X_f)$ denote its coordinate ring. 
\begin{enumerate}
\item The {\em semistable locus} for the $U$-action on $X$ linearised by $L$ is $X^{ss,
U} =  \bigcup_{f \in I^{\mathrm{fg}}} X_f$ where
$$I^{\mathrm{fg}} = \{f
\in I \ | \ {\cO}(X_f)^U
\mbox{ is finitely generated} \}.$$
\item The {\em (locally trivial) stable  locus} for the linearised $U$-action on $X$ is $ X^{\mathrm{lts},U} =
\bigcup_{f \in I^{\mathrm{lts}} } X_f$ where 
\[ I^{\mathrm{lts}}:= \{f
\in I^{\mathrm{fg}}  \ | \ 
 q_U: X_f \ra \Spec({\cO}(X_f)^U) \mbox{ is a locally trivial
geometric quotient} \}.\]
\item \label{defn:envelopquot}
The {\em enveloped quotient} of
$X^{ss,U}$ by the linear $U$-action is $q_U: X^{ss, U} \rightarrow q_U(X^{ss,U})$, where
$q_U: X^{ss, U} \rightarrow \Proj({\widehat{\cO}}_L(X)^U)$ is the natural
morphism of schemes and
$q_U(X^{ss,U})$ is a dense constructible subset of the {\em
enveloping quotient}
$$X \env  U = \bigcup_{f \in I^{fg}}
\Spec({\calo}(X_f)^U).$$ 
\end{enumerate}
\end{definition}

\begin{rmk}\label{remclaim}\
\begin{enumerate}
\item Even when ${\widehat{\calo}}_L(X)^U$ is finitely generated, so that  $X\env U=\Proj({\widehat{\calo}}_L(X)^U)$, the enveloped quotient $q_U(X^{ss,U})$ is not necessarily a subscheme of $X \env U $ (for example, see \cite[$\S$6]{DK}).
\item The enveloping quotient 
$X \env  U$ has quasi-projective open subschemes (\lq inner enveloping quotients' $X \inenv U$) that contain the enveloped quotient $q_U(X^{ss})$  and have ample line bundles which under the natural map $q_U:X^{ss} \to X\env U$ pull back to positive tensor powers of $L$ (see \cite{BDHK} for details). 
\end{enumerate}
\end{rmk}

The $H$-semistable locus $X^{ss} = X^{ss,H}$ and enveloped and (inner) enveloping quotients
$$q_H: X^{ss} \to q_H(X^{ss}) \subseteq X \inenv H \subseteq X \env H$$
for the linear action of $H$ are defined as for the unipotent case in Definition \ref{defnssetc} and Remark \ref{remclaim} (cf. \cite{BDHK}), but the definition of the stable locus $X^{\mathrm{lts}}= X^{\mathrm{lts},H}$ for the linear action of $H$ combines (and extends) the definitions for unipotent and reductive groups.

\begin{definition} (Stability for linear algebraic groups). \label{def:GiSt1.1}
Let $H$ be a linear algebraic group acting linearly on $X$ with respect to an ample line bundle $L$; then the \emph{(locally trivial) stable locus} is the open subscheme $X^{\mathrm{lts}}= \bigcup_{f \in I^{\mathrm{lts}}} X_f$ of $X^{ss}$, where $I^{\mathrm{lts}} \subseteq \bigcup_{r>0} H^0(X,L^{\ten r})^H$ is the subset
of $H$-invariant sections $f$ satisfying the following conditions:
\begin{enumerate}
\item \label{itm:GiSt1.1-1} the $H$-invariant open subscheme $X_f$ is affine;
\item \label{itm:GiSt1.1-2} the $H$-action on $X_f$ is closed with finite stabiliser groups; and
\item \label{itm:GiSt1.1-3} the restriction of the $U$-enveloping quotient map 
 \[
q_{U}:X_f \to \spec(({\widehat{\calo}}_L(X)^{U})_{(f)})
\]
is a locally trivial geometric quotient for the $U$-action on $X_f$. 
\end{enumerate}
\end{definition}

\subsection{GIT for linear algebraic groups with internally graded unipotent radicals}\label{GIT_int_gr}

Now suppose that $H = U \rtimes R$ is a linear algebraic group with internally graded unipotent radical $U$ in the following sense.

\begin{definition}\label{int gr uni rad} We say $H = U \rtimes R$ has {\em internally graded unipotent radical} $U$ if there is a central 1-PS $\lambda: \GG_m \to R$ whose conjugation action on the Lie algebra of $U$ has strictly positive weights. We write $\hU = U \rtimes \lambda(\GG_m) \leqslant H$ for the associated semi-direct product.
\end{definition}

Suppose also that ${H}$ acts linearly on $X$ with respect to an ample line bundle $L$.  It is shown in \cite{BDHK2, BDHK3} that the algebra of $H$-invariant sections is finitely generated provided:
\begin{enumerate}[label={\alph*})]
\item $L$ is replaced with a suitable tensor power $L^{\otimes m}$, with $m\geq 1$ sufficiently divisible, and  the linearisation of the action of ${{H}}$ is twisted by a suitable (rational) character, and 
\item condition \eqref{star} described below (also known as \lq semistability coincides with stability for the unipotent radical $U$') holds.
\end{enumerate}
\noindent  Moreover, in this situation the natural quotient morphism $q_H$ from the semistable locus $X^{ss,{{H}}}$ to the enveloping quotient $X\env {{H}}$ is surjective, and expresses the projective scheme $X/\!/H = X \inenv H = X\env {{H}}$ as a good quotient of $X^{ss,{{H}}}$. Furthermore this locus $X^{ss,{{H}}}$  can be described using Hilbert--Mumford criteria. It is also shown in \cite{BDHK3} that when condition \eqref{star} is not satisfied, but is replaced with a slightly weaker condition, such as \eqref{star2} below, then there is a sequence of blow-ups of $X$ along $H$-invariant subschemes (similar to that of \cite{K2} when $H$ is reductive) resulting in a projective scheme $\widehat{X}$ with an induced linear action of ${H}$ satisfying condition \eqref{star}. In fact, these results can be generalised to allow actions where the $U$-action has positive dimensional stabilisers (\textit{cf.}\ Theorems \ref{Uh3} and \ref{Uh4} below). Before giving a precise description of the condition \eqref{star} and its variants, we define the notion of an adapted linearisation, which is also needed for the statement of the main results of \cite{BDHK2, BDHK3}. 

Let $\chi: H  \to \GG_m$ be a character of $H$; the restriction to $\hU$ of $\chi$ contains $U$ in its kernel and can be identified naturally with an integer so that the integer 1 corresponds to the character of $\hU$ which fits into the exact sequence $U \hookrightarrow \hU \twoheadrightarrow \lambda(\GG_m)$. By replacing $L$ with a sufficiently high power, we can without loss of generality assume that $L$ is very ample. Let $\omega_{\min}  := \omega_0 < \omega_{1} < \cdots < \omega_{\max} $ be the weights with which the 1-PS $ \lambda : \GG_m \ra {\hU} \leq \widehat{H}$ acts on the fibres of the tautological line bundle $\mathcal{O}_{\PP((H^0(X,L)^*)}(-1)$ over points of the fixed locus $\PP((H^0(X,L)^*)^{\lambda(\GG_m)}$. Without loss of generality we may assume that there exist at least two distinct such weights, as otherwise the $U$-action on $X$ is trivial, in which case we can take a quotient by the action of the reductive group $R=H/U$. 

\begin{definition}\label{adapted}
For a character $\chi$ of $H$ as above and a positive integer $c$, we say the rational character $\chi/c$ is {\em  adapted} to the linear action of $\hU$ if
\begin{equation} 
\omega_{\min}:=\omega_{0} < \frac{\chi}{c} < \omega_{1}. 
\end{equation}
Furthermore, we say $L$ is adapted to the $\hU$-action if $\omega_{\min}:=\omega_0 <0 < \omega_1 $.
\end{definition}

If the rational character $\chi/c$ is adapted to the linear action of $\hU$, then the $H$-linearisation $\mathcal{L}_\chi^{\otimes c}$ on $X$ given by twisting the ample line bundle $L^{\otimes c}$ by the character $\chi$ (that is, so that the weights $\omega_j$ are replaced with $\omega_jc-\chi$) is adapted. Let $X^{s,\GG_m}_{\min+}=X^{s,\lambda(\GG_m)}_{\min+}$ denote the stable set in $X$ for the linear action of $\GG_m$ via $\lambda$ with respect to the adapted linearisation $\mathcal{L}_\chi^{\otimes c}$ and, for a maximal torus $T$ of $H$ containing $\lambda(\GG_m)$, let $X^{s,T}_{\min+}$ denote the stable set in $X$ for the linear action of $T$ with respect to the adapted linearisation $L_\chi^{\otimes c}$. By the theory of variation of (classical) GIT \cite{dh98,Thaddeus}, the stable set $X^{s,\lambda(\GG_m)}_{\min+}$ is independent of the choice of adapted rational character $\chi/c$. In fact, by the Hilbert--Mumford criterion, we have $X^{s,\lambda(\GG_m)}_{\min+}  = X^0_{\min} \setminus Z_{\min}$, where if $V_{\min}$ denotes the weight space  of $\omega_{\min}$ in $V=H^0(X,L)^*$, then
\[
Z_{\min} :=X \cap \PP(V_{\min})=\left\{x \in X^{\lambda(\GG_m)} : \lambda(\GG_m) \text{ acts on } L^*|_x \text{ with weight } \omega_{\min} \right\} \]
and
\[
X^0_{\min} :=\{x\in X \mid \lim_{t \to 0, \,\, t \in \GG_m} \lambda(t) \cdot x \in Z_{\min}\}.
\]   

\begin{definition} \label{cond star}(Conditions $(*) - (***)$ generalising \lq semistability equals stability' \textit{cf.}\ \cite{BDHK3}) 
With the above notation, we define the following conditions for the $\hU$-action on $X$.
\begin{equation}\label{star}
 \stab_{U}(z) = \{ e \}  \textrm{ for every } z \in Z_{\min}. \tag{$*$}
\end{equation}
\begin{equation}\label{star2}
 \stab_{U}(x) = \{ e \}  \textrm{ for generic } x \in X_{\min}^0. \tag{$**$}
\end{equation}
Moreover, if $ U \geqslant U^{(1)} \geqslant ... \geqslant U^{(s)}  \geqslant \{e\}$ denotes the derived series of $U$, we define condition
\begin{equation}\label{star3}
\textrm{for } 1 \leq j \leq s, \textrm{ there exists } d_j \in \mathbb{N} \textrm{ such that } \dim \stab_{U^{(j)}}(x) = d_j \textrm{ for all } x \in X^0_{\min}. \tag{$***$}
\end{equation}
Note that condition \eqref{star} holds if and only if we have $\stab_{U}(x) = \{e\}$ for all $x \in X^0_{\min}$. This condition is also referred to  in \cite{BDHK3} as the condition \lq semistability coincides with stability' for the action of $\hU$ (or, when the 1-PS $\lambda:\GG_m \to R$ is fixed, for the linear action of $U$). 
\end{definition}

\begin{definition}\label{def min stable}
Let $T \leq R$ be a maximal torus containing $\lambda(\GG_m)$. The {\em minimal stable set} for a linear $H$-action for which condition \eqref{star} holds for the action of the graded unipotent group $\hU$ is
\[X^{s,{{H}}}_{\min+} :=  \bigcap_{h \in H} h X^{s,T}_{\min+}.\]
\end{definition}

We note that the minimal stable set for the $\hU$-action satisfies 
\[ X^{s,{\hU}}_{\min+}=  \bigcap_{u \in U} u X^{s,\lambda(\GG_m)}_{\min+} = X^0_{\min} \setminus U Z_{\min}.  \]

\begin{thm} \label{mainthm} 
 \textbf{($\hU$-Theorem when semistability coincides with stability for $\hU$)}  {\rm \cite{BDHK3}}
Suppose that the linearisation for the action of $H$ on $X$ is adapted as in Definition \ref{adapted} and that the $\hU$-action on $X$ satisfies condition \eqref{star}. Then the following statements hold.
\begin{enumerate}[label=(\emph{\roman*})]
\item The open subscheme $X^{s,{\hU}}_{\min+}$ of $X$ has a projective geometric quotient $X /\!/ \hU = X^{s,{\hU}}_{\min+}/\hU$ by $\hU$.
\item The open subscheme $X^{s,{{H}}}_{\min+}$ of $X$ has a good quotient $X /\!/ H = (X /\!/ \hU)/\!/R$ by $H=U \rtimes R$, which is also projective.
\end{enumerate}  
\end{thm} 

\begin{rmk} In the proof of Theorem \ref{mainthm} (and its variants below), one replaces the adapted linearisation by a \lq well adapted' linearisation (which can be achieved by twisting by a rational character); this is a slightly stronger notion. This strengthening does not alter $X^{s,{\hU}}_{\min+}$ or its quotient $X/\!/ \hU$, but it affects what can be said about induced ample line bundles on $X /\!/  \hU$ and $X/\!/ H = (X/\!/ \hU)/\!/ (R/\lambda(\GG_m))$.
The proofs in \cite{BDHK2,BDHK3} that the algebras of invariants $\oplus_{m \geq 0} H^0(X,L_{m\chi}^{\otimes cm})^{\hU}$ and 
\[\bigoplus_{m\geq 0} H^0(X,L_{m\chi}^{\otimes cm})^{{H}} = (\bigoplus_{m \geq 0} H^0(X,L_{m\chi}^{\otimes cm})^{\hU})^{(R/\lambda(\GG_m))}\]
are finitely generated, and that the enveloping quotients  $X /\!/  \hU = X^{s,{\hU}}_{\min+}/\hU$ and $X /\!/ H$ are the associated projective schemes, require that the linearisation is twisted by a  well adapted rational character $\chi/c$. More precisely, it is shown in \cite{BDHK2, BDHK3} that, given a linear action of $H$ on $X$ with respect to an ample line bundle $L$, there exists $\epsilon >0$ such that if $\chi/c$ is a rational character of $\GG_m$ (lifting to $H$) with $c$ sufficiently divisible and $\chi: {{H}} \to \GG_m$ a character of $H$ such that
$$\omega_{\min} < \frac{\chi}{c} < \omega_{\min} + \epsilon,$$
then the algebras of invariants $\oplus_{m\geq 0} H^0(X,L_{m\chi}^{\otimes cm})^{\hU}$ and $\oplus_{m\geq 0} H^0(X,L_{m\chi}^{\otimes cm})^{{H}} $ are finitely generated, and the associated projective schemes 
$X /\!/  \hU$ and $X /\!/ H$ satisfy the conclusions of Theorem \ref{mainthm}.
\end{rmk}

Theorem \ref{mainthm} describes the good case when semistability coincides with stability for the linear action of $\hU$. The following versions proved in \cite{BDHK3} apply more generally.

\begin{thm}  \label{Uh2}
 \textbf{($\hU$-Theorem giving projective completions)}  {\rm \cite{BDHK3}}
Suppose that the linear action of ${\hU}$ on $X$ is adapted and satisfies condition \eqref{star2}. Then there exists a sequence of blow-ups along $H$-invariant projective subschemes, resulting in a projective scheme $\widetilde{X}$ (with blow-down map $\widetilde{\psi}: \widetilde{X} \to X$)  such that the conditions of Theorem \ref{mainthm} still hold for a suitable ample linearisation, and such that the quotient  given by that theorem is a geometric quotient of an open subscheme $\widetilde{X}^{s,H}$ of $\widetilde{X}$.
\end{thm}

\begin{rmk} \label{centbu}
This blow-up sequence is constructed in two stages: one first blows up by considering the stabiliser subgroups for the unipotent group $U$, and then one blows up by considering the stabiliser subgroups for the reductive group $R$. 

In the first step, one performs a blow-up sequence to obtain $\widehat{\psi}: \widehat{X} \to X$ such that the $\hU$-action on $\widehat{X}$ satisfies condition \eqref{star} with respect to a linearisation $\widehat{\cL}$, which is an arbitrarily small perturbation of $\widehat{\psi}^*(\mathcal{L})$. The centres of the blow-ups used to obtain $\widehat{X}$ from $X$ are determined by the dimensions of the stabilisers in $U$ of the limits $\lim_{t \to 0} \lambda(t) \cdot x$ for $x \in X^{0}_{\min}$ (for details, see \cite{BDHK3}). Then one can construct a projective and geometric quotient of the $\hU$-action on $\widehat{X}^{s,\hU}$ by Theorem \ref{mainthm} and, as $\widehat{\psi}$ is an isomorphism away from the exceptional divisor, one obtains a geometric quotient of a $\hU$-invariant open subset $X^{{s},\hU}$ of $X$ as an open subscheme of $\widehat{X}  /\!/ _{\widehat{\mathcal{L}}} \hU$, where $X^{{s},\hU}$ is the image under $\widehat{\psi}$ of the intersection of $\widehat{X}^{s,\hU}$ with the complement of the exceptional divisor in $\widehat{X}$. Another characterisation of this stable locus is as
$$X^{s,\hU} = \{ x \in \psi (\widehat{X}^{s,\hU})  \mid \dim \stab_U(\lim_{t \to 0} \lambda(t) \cdot x) = 0 \}. $$
If one is only interested in obtaining a good quotient for the $H$-action, then the second stage of the blow-up procedure is not needed: one can then take a reductive GIT quotient of the residual action on $\widehat{X}  /\!/\!_{\widehat{\mathcal{L}}} \hU$ of $H/\hU = R/\lambda(\GG_m)$, and thus one obtains a good quotient of the $H$-action on an open subset of $X$ as an open subscheme of $(\widehat{X}  /\!/ _{\widehat{\mathcal{L}}} \hU)/\!/ (R/\lambda(\GG_m))$. Moreover, this good quotient restricts to a geometric quotient on an open subscheme of stable points.

Finally, to go from $\widehat{X}$ to the blow-up $\widetilde{X}$ in Theorem \ref{Uh2}, one performs an additional blow up sequence by considering the stabiliser groups for the action of the reductive group $R/\lambda(\GG_m)$ as in the partial desingularisation procedure described in $\S$\ref{sec part desing}.

This gives us an $H$-invariant open subscheme $X^{{s},H}$ of $X$ with a geometric quotient by $H$ which is an open subscheme of the projective scheme $\widetilde{X}  /\!/\!_{\widetilde{\mathcal{L}}} H$ (and also of $\widehat{X}  /\!/\!_{\widehat{\mathcal{L}}} H$). Here $X^{{s},H}$ is the image under $\widehat{\psi}$ of the intersection of $\widehat{X}^{s,H}$ with the complement of the exceptional divisor in $\widehat{X}$.
\end{rmk}

Theorem \ref{mainthm} can be generalised by weaking the condition \eqref{star} further to \eqref{star3} to allow for actions with positive dimensional stabiliser groups generically.

\begin{thm} \textbf{($\hU$-Theorem with positive-dimensional stabilisers in $U$)}  {\rm \cite{BDHK3}}
\label{Uh3}
Suppose that condition \eqref{star3} holds for an adapted linear $\hU$-action on $X$. Then the conclusions of Theorem \ref{mainthm} hold. 
\end{thm}

In fact, this theorem still holds if we replace the derived series in condition \eqref{star3} with any series $U \geqslant U^{(1)} \geqslant ... \geqslant U^{(s)}  \geqslant \{e\}$ which is normalised by $H$ and whose successive quotients $U^{(j)}/U^{(j+1)}$ are abelian, provided that \eqref{star3} holds for this series.  

Finally, there is a version of the theorem without requiring any hypothesis related to semistability coinciding with stability.

\begin{thm} \label{Uh4} \textbf{($\hU$-Theorem with positive-dimensional stabilisers in $U$, giving projective completions)}  {\rm \cite{BDHK3}}
For a linear $H$-action on $X$ with respect to an adapted ample linearisation $\cL$, there is a sequence of blow-ups along $H$-invariant projective subschemes, resulting in a projective scheme $\widetilde{X}$ such that condition \eqref{star3} holds for a suitable linearisation, and such that the $H$-quotient given by Theorem \ref{mainthm} is a geometric quotient of an open subscheme $\widetilde{X}^{s,H}$ of $\widetilde{X}$.
\end{thm}

This theorem provides a non-reductive analogue of the partial desingularisation construction for reductive GIT described at the end of $\S$\ref{sec GIT classical}. 

As before, this gives us an $H$-invariant open subscheme $X^{s,H}$ of $X$ with a geometric quotient by $H$ which is open in the projective scheme $X /\!/\!_{\mathcal{L}} H$; here $X^{s,H}$ is the image under $\psi$ of the intersection of $\widetilde{X}^{s,H}$ with the complement of the exceptional divisor in $\widetilde{X}$.

\begin{rmk} \label{centbu2} 
As for Theorem \ref{Uh2} (\textit{cf.}\ Remark \ref{centbu}), one first constructs a blow-up $\widehat{X} \ra X$ by considering stabiliser subgroups for the unipotent action, then one constructs a further blow-up sequence $\widetilde{X} \ra \widehat{X}$ by considering the stabiliser subgroups for the residual action of the reductive subgroup $R/\lambda(\GG_m)$. The slight difference is that in the first step, the centres of the blow-ups used to obtain $\widehat{X}$ from $X$ are determined by the dimensions of the $U^{(j)}$-stabilisers of the limit $\lim_{t \to 0} \lambda(t) \cdot x$ and of $x$ itself for $x \in X^{0}_{\min}$. We will see in the next section that by applying Theorem \ref{Uh4} to the closures of the subschemes where these dimensions take different values, and combining this with the partial desingularisation construction of \cite{K2} for reductive GIT quotients,  $X$ can be stratified so that each stratum is a locally closed $H$-invariant subscheme of $X$ with a  geometric quotient by the action of $H$.
\end{rmk}

\section{Stratifying quotient stacks}\label{sec strat quot stacks}

Let $H = U \rtimes R$ be a linear algebraic group with internally graded unipotent radical $U$ (\textit{cf.}\ Definition \ref{int gr uni rad}) acting on a projective scheme $X$ over an algebraically closed field $\kk$ of characteristic 0 with respect to an ample linearisation $\cL$. We fix an invariant inner product on the Lie algebra $\Lie R$ of the Levi factor, just as in the construction of the HKKN stratification (\textit{cf.}\ $\S$\ref{stratification}).

The aim of this section is to prove Theorem \ref{mainthmstrat} stated in $\S$\ref{sec intro}. We will prove this result using a recursive argument involving the dimensions of $X$ and of $H$ and the number of irreducible components of $X$. The idea will be to start by defining a \lq minimum' stratum, which will be a non-empty $H$-invariant open subscheme of $X$, and then proceed recursively. 

We will assume that $H$ is connected and that $X$ is reduced.

\begin{rmk} These assumptions do not involve any significant loss of generality.\
\begin{enumerate}
\item
The HKKN stratification $\{S_\beta \mid \beta \in \mathcal{B} \}$ is usually indexed by a finite subset $\mathcal{B}$ of a positive Weyl chamber. Then the strata are not necessarily connected even when $X$ is irreducible, and it is often useful to refine the stratification so that the strata are the connected components of $S_\beta$, or are unions of some but not all of these connected components (\textit{cf.}\ \cite{K}). There is a similar ambiguity in the construction of the refined stratification defined in this section: at some points we take connected components, but this is not crucial to the definition. Indeed if we wish to allow the group $H$ to be disconnected then we cannot assume that the strata are connected since they are required to be $H$-invariant. Then instead of taking connected components (which will be invariant under the component $H_0$ of the identity in $H$), we can take their $H$-sweeps, which will be disjoint unions of at most $|H/H_0|$ of these connected components.
\item
If $X$ is non-reduced, then we can define the stratification on $X$ by using a positive power of $L$ to define a $H$-equivariant embedding of $X$ into a projective space  $\PP^n$, and then take the fibre product of $X$ with the stratification on $\PP^n$. Indeed, we will see that the stratification is functorial for equivariant closed immersions (this is essentially the third statement in Theorem \ref{mainthmstrat}). This follows as the reductive notions of GIT (semi)stability are functorial and since in the non-reductive case, we are assuming that $H$ has internally graded unipotent radical. The stable loci when $H$ has internally graded unipotent radical and adapted linearisation are also functorial, as they have Hilbert--Mumford style descriptions (see Definition \ref{def min stable}). We note that in the more general non-reductive GIT set up described in $\S$\ref{sec non-red} the notions of (semi)stability are not functorial, as taking $H$-invariants is not exact, and so there can be invariants which do not extend to the ambient space. The advantage of assuming that $X$ is reduced is that the complement to the open stratum then has a canonical scheme structure; thus it is easier to recursively define the stratification.
\end{enumerate}
\end{rmk}

Let us describe the recursive construction when $H$ is connected and $X$ is reduced. For each linear action on $X$ of $H= U \rtimes R$ with internal grading $\lambda: \GG_m \to R$ and linearisation $\mathcal{L}$ with underlying ample line bundle $L$, we will first use recursion to define a {\em nonempty} $H$-invariant open subscheme $\mathcal{S}_0(X,H,\lambda,\mathcal{L})$ of $X$ that admits a geometric quotient $\mathcal{S}_0(X,H,\lambda,\mathcal{L})/H$; this will be done by considering seven different cases. After defining the open stratum, we will define the stratification $\{ \mathcal{S}_\gamma | \gamma \in \Gamma \}$ of $X$ with strata $\mathcal{S}_\gamma = \mathcal{S}_\gamma(X,H,\lambda,\mathcal{L})$ and index set $\Gamma = \Gamma(X,H,\lambda,\mathcal{L})$ by letting $X_1, \ldots,X_k$ be the connected components of the projective subscheme of $X$ equal to the complement of $\mathcal{S}_0(X,H,\lambda,\mathcal{L})$, letting
$\mathcal{S}_{0,i}(X,H,\lambda,\mathcal{L})$ for $1 \leqslant i \leqslant m$ be the connected components of $\mathcal{S}_0(X,H,\lambda,\mathcal{L})$,
and setting
\begin{equation} \label{(1)} \Gamma(X,H,\lambda,\mathcal{L}) := \{0\}\times \{1,\ldots,m\}  \cup \bigcup_{1 \leqslant j \leqslant k} \{ X_j\} \times \Gamma(X_j ,H,\lambda,\mathcal{L}|_{X_j}).
\end{equation}
The strata indexed by $(0,i) \in \Gamma(X,H,\lambda,\mathcal{L})$ for $1 \leqslant i \leqslant m$ are then the connected components of the open subscheme  the open subscheme $\mathcal{S}_0(X,H,\lambda,\mathcal{L})$ constructed using the case by case argument below, whereas the stratum indexed by an element $(X_j, \gamma)$ for $1 \leq j \leq k$ and $\gamma \in
\Gamma(X_j,H,\lambda,\mathcal{L}|_{X_j})$ is 
\begin{equation} \label{(2)}  \mathcal{S}_{(X_j,\gamma)}(X,H,\lambda,\mathcal{L}) := \mathcal{S}_\gamma(X_j,H,\lambda,\mathcal{L}|_{X_j}),
\end{equation}
where the strata $\mathcal{S}_\gamma(X_j,H,\lambda,\mathcal{L}|_{X_j})$ are constructed by induction. The partial order on $\Gamma$ then naturally comes from the partial orders on each $\Gamma(X_j ,H,\lambda,\mathcal{L}|_{X_j})$ 
with $(0,i) < (X_j, \gamma)$ for all $1 \leqslant i \leqslant m$ and $1 \leqslant j \leqslant k$ and $\gamma \in \Gamma(X_j ,H,\lambda,\mathcal{L}|_{X_j})$.

Let us now describe how to define $\mathcal{S}_0(X,H,\lambda,\mathcal{L})$ in the seven different cases.
Let  $$ U = U^{(0)} \geqslant U^{(1)} = [U,U] \geqslant ... \geqslant U^{(s)}  \geqslant \{e\}$$ be the derived series of $U$. Let
$$Z_{\min}^{d_0} = \{ x \in Z_{\min} | \dim(\text{Stab}_{U}(x)) = d_0 \}$$
where $d_0$ is the minimal value of $ \dim(\text{Stab}_{U}(x))$ for $x \in Z_{\min}$. Then $Z_{\min}^{d_0}$ is a nonempty $H$-invariant open subscheme of $Z_{\min}$, and $UZ_{\min}^{d_0}$ is a nonempty locally closed subscheme of $X$. 

\medskip

\noindent {\bf Case 1.}
First assume that the central one-parameter subgroup $\lambda: \GG_m \to R$ of $R$ has at least two distinct weights for the linear action of $H$ on $X$; equivalently $Z_{\min}$ is a projective subscheme of $X$ with $Z_{\min} \neq X$. Under this assumption we have two possibilities to consider.

\medskip

\noindent {\bf Case 1(a).}  Suppose that $UZ_{\min}^{d_0}$ is open in $X$.
Then we define $$\mathcal{S}_0(X,H,\lambda,\mathcal{L}) = U \{x \in \mathcal{S}_0(Z_{\min},R/\lambda(\GG_m),\lambda_0,\mathcal{L}|_{Z_{\min}})  \,\, | \,\, \dim(\text{Stab}_{U_j}(x)) = d_j \mbox{ for } 1 \leqslant j \leqslant s \}
$$
where $\lambda_0$ is the trivial one-parameter subgroup that grades the trivial unipotent radical of the reductive group $R/\lambda(\GG_m)$ and $d_j:= \min\{\dim(\text{Stab}_{U_j}(x)) :x \in \mathcal{S}_0(Z_{\min},R/\lambda(\GG_m),\lambda_0,\mathcal{L}|_{Z_{\min}}) \}$. 

By induction we can assume that $ \mathcal{S}_0(Z_{\min},R/\lambda(\GG_m),\lambda_0,\mathcal{L}|_{Z_{\min}}) $ is a nonempty $R$-invariant open subscheme of $Z_{\min}$ with a geometric quotient by the action of $R/\lambda(\GG_m)$ (or equivalently by the action of $R$, since the central one-parameter subgroup $\lambda(\GG_m)$ of $R$ acts trivially on $Z_{\min}$). Indeed, we can construct such an open subscheme as in Case 2 described below. Thus $\mathcal{S}_0(X,H,\lambda,\mathcal{L})$ is an $H$-invariant nonempty open subscheme of $X$, and by the proof of Theorem \ref{Uh4} (see \cite[Remark 2.10]{BDHK3}) 
it has a geometric quotient
\begin{align*}
\mathcal{S}_0(X,H,\lambda,\mathcal{L})/H & \cong \mathcal{S}_0(Z_{\min},R/\lambda(\GG_m),\lambda_0,\mathcal{L}|_{Z_{\min}}) /R \\
 & = \mathcal{S}_0(Z_{\min},R/\lambda(\GG_m),\lambda_0,\mathcal{L}|_{Z_{\min}})  /(R/\lambda(\GG_m)).
\end{align*}

\medskip

\noindent {\bf Case 1(b).} Suppose that $UZ_{\min}^{d_0}$ is not open in $X$. Recall that $X^0_{\min}$ is an $H$-invariant open subscheme of $X$ and that the morphism $p:X^0_{\min} \to Z_{\min}$ defined by
$$ p(x) = \lim_{t \to 0} \lambda(t) x$$
satisfies $p(urx) = rp(x)$ for all $x \in X^0_{\min}$, $u \in U$ and $r \in R$ (\textit{cf.}\ \cite{BDHK3}). We set 
\[\mathcal{S}_0(X,H,\lambda,\mathcal{L})= \left\{x \in X^0_{\min} \setminus U Z_{\min} \: | \:  \begin{array}{c}
p(x) \in \mathcal{S}_0(Z_{\min},R/\lambda(\GG_m),\lambda_0,\mathcal{L}|_{Z_{\min}}) \mathrm{ \: and \: for \:} 1 \leqslant j \leqslant s \\ \dim(\text{Stab}_{U_j}(x)) = \delta_j  \mbox{ and } \dim(\text{Stab}_{U_j}(p(x))) = d_j
\end{array} \right\}\]
where $d_j$ and $\delta_j$ are defined inductively, in decreasing order of $j$, as follows. Assuming that $d_i$ and $\delta_i$ are defined for $i > j$, we define 
\[ d_j := \min\left\{ \dim(\text{Stab}_{U_j}(z)) \: | \: \begin{array}{c} \exists \: x \in X^0_{\min} \textrm{ with }z= p(x) \in \mathcal{S}_0(Z_{\min},R/\lambda(\GG_m),\lambda_0,\mathcal{L}|_{Z_{\min}}) \\ \textrm{and } \dim(\text{Stab}_{U_i}(z)) = d_i \textrm{ for all } i > j \end{array} \right\} \]
and
\[ \delta_j := \min\left\{ \dim(\text{Stab}_{U_j}(x)) \: | \: \begin{array}{c}  x \in X^0_{\min} \textrm{ such that } p(x) \in \mathcal{S}_0(Z_{\min},R/\lambda(\GG_m),\lambda_0,\mathcal{L}|_{Z_{\min}}) \\ \textrm{and } \dim(\text{Stab}_{U_i}(x)) = \delta_i \textrm{ for all } i > j \end{array} \right\}, \]
where for $j =s$, the conditions appearing in the second line of these definitions are empty.

Note that if $X$ and $Z_{\min}$ are irreducible then we can define $d_j$ and $\delta_j$ much more simply as the generic values taken by $\dim(\text{Stab}_{U_j}(z)) $ for $z \in Z_{\min}$ and by $\dim(\text{Stab}_{U_j}(x)) $ for $x \in X^0_{\min}$, but in this approach to defining the stratification  $\{ \mathcal{S}_\gamma | \gamma \in \Gamma \}$ we need to allow $X$ to be reducible.

Using the proof of Theorem \ref{Uh4} again (cf. \cite[Remark 2.10]{BDHK3}), we have that $\mathcal{S}_0(X,H,\lambda,\mathcal{L})$ is an $H$-invariant nonempty open subscheme of $X$ with a geometric quotient
$\mathcal{S}_0(X,H,\lambda,\mathcal{L})/\hU$, and by the inductive construction a geometric quotient $$\mathcal{S}_0(X,H,\lambda,\mathcal{L})/H = (\mathcal{S}_0(X,H,\lambda,\mathcal{L})/\hU)/(R/\lambda(\GG_m)).$$

\medskip

\noindent {\bf Case 2.} Now assume that the central one-parameter subgroup $\lambda: \GG_m \to R$ which grades $U$ acts trivially on $X$, so the unipotent radical $U$ of $H = U \rtimes R$ must also act trivially on $X$. Then without loss of generality we can assume that $H=R$ is reductive and that $\lambda = \lambda_0$ is trivial.

\medskip

\noindent {\bf Case 2(a).} Suppose that the stable locus $X^{s,R}$ for the linear action of $R$ on $X$ is nonempty. Then we let 
$$\mathcal{S}_0(X,H,\lambda,\mathcal{L}) = X^{s,R}$$
and by classical GIT this has a geometric quotient $X^{s,R}/H = X^{s,R}/R$.

\medskip

\noindent {\bf Case 2(b).} Suppose that the semistable locus $X^{ss,R}$ for the linear action of $R$ on $X$ is empty. Then there is a stratum $S_\beta$ from the HKKN stratification for $X$ (associated to our invariant inner product on $R$) such that $\beta \neq 0$ and
$$S_\beta = HY^{ss}_\beta = RY^{ss}_\beta \cong R \times^{P_\beta} Y^{ss}_\beta$$
is nonempty and open in $X$
(see $\S$\ref{stratification} and Remark \ref{rem:specialcases}).  Then we have the following two subcases to consider.

\medskip

\noindent {\bf Case 2(b)i).} Suppose that $\overline{Y^{ss}_\beta} \neq X$. Then by induction on $\dim X$ and the number of irreducible components of $X$, 
$$\mathcal{S}_0(X,H,\lambda,\mathcal{L}) = R( \mathcal{S}_0(\overline{Y^{ss}_\beta},P_\beta,\lambda_\beta,\mathcal{L}|_{\overline{Y^{ss}_\beta}}) \cap Y^{ss}_\beta) \cong R \times^{P_\beta} ( \mathcal{S}_0(\overline{Y^{ss}_\beta},P_\beta,\lambda_\beta,\mathcal{L}|_{\overline{Y^{ss}_\beta}}) \cap Y^{ss}_\beta)
$$
is a nonempty $R$-invariant (and hence $H$-invariant) open subscheme of $X$ and has a geometric quotient 
$$\mathcal{S}_0(X,H,\lambda,\mathcal{L})/H = \mathcal{S}_0(X,H,\lambda,\mathcal{L})/R  \cong  ( \mathcal{S}_0(\overline{Y^{ss}_\beta},P_\beta,\lambda_\beta,\mathcal{L}|_{\overline{Y^{ss}_\beta}}) \cap Y^{ss}_\beta)/P_\beta.
$$

\medskip

\noindent {\bf Case 2(b)ii).} Suppose that $\overline{Y^{ss}_\beta} = X$. Then $P_\beta = R$ so $\beta$ defines a rational character of $R$ and corresponds to a nontrivial central one-parameter subgroup $\lambda_\beta: \GG_m \to R$. If $Z_\beta \neq X$ then the one-parameter subgroup $\lambda_\beta(\GG_m)$ acts nontrivially on $X$, and thus we can use Case 1 and induction to define
$$ \mathcal{S}_0(X,H,\lambda,\mathcal{L})  = \mathcal{S}_0(X,R,\lambda_\beta,\mathcal{L}) $$
so that it is a nonempty $R$-invariant (and hence $H$-invariant) open subscheme of $X$ and has a geometric quotient
$$ \mathcal{S}_0(X,H,\lambda,\mathcal{L}) / H = \mathcal{S}_0(X,R,\lambda_\beta,\mathcal{L}) /R.$$
If $Z_\beta = X$ then $\lambda_\beta(\GG_m)$ acts trivially on $X$ and we can use induction on the dimension of $H$
to define
$$\mathcal{S}_0(X,H,\lambda,\mathcal{L})  = \mathcal{S}_0(X,R/\lambda_\beta(\GG_m),\lambda_0,\mathcal{L}) $$
which is a nonempty $R$-invariant (and hence $H$-invariant) open subscheme of $X$ with geometric quotient
$ \mathcal{S}_0(X,H,\lambda,\mathcal{L}) / H = \mathcal{S}_0(X,R/\lambda_\beta(\GG_m),\lambda_0,\mathcal{L}) /(R/\lambda_\beta(\GG_m)).$

\medskip

\noindent {\bf Case 2(c).} Suppose now that $X^{s,R} = \emptyset \neq X^{ss,R}$. Recall from Remark \ref{rem:specialcases} that the partial desingularisation construction \cite{K2} for the linear action of the reductive group $R$ on $X$ can be applied to $X^{ss,R}$, although, for simplicity, $X$ was assumed to be irreducible in Remark \ref{rem:specialcases}, which is no longer the case here.  This construction terminates with a birational projective morphism $\psi: \widetilde{X} \to X$ and $R$-linearisation $\widetilde{\mathcal{L}}$ for an ample line bundle $\widetilde{L}$ on $\widetilde{X}$ which is a positive integer multiple of $\psi^*L \otimes \mathcal{O}(-\epsilon E)$ where $E$ is the exceptional divisor and $0< \epsilon <\!< 1$ is rational. Then we have the following two subcases to consider.

\medskip

\noindent {\bf Case 2(c)i).} Suppose that 
 there is a positive-dimensional connected reductive subgroup $R'$ of $R$ such that 
$$R \widetilde{Z}_{R'}^{ss} \cong  R \times^{N_{R'}} \widetilde{Z}_{R'}^{ss}$$
is open and nonempty in $\widetilde{X}^{ss,R}$. Here $N_{R'}$ is the normaliser of $R'$ in $R$ and 
$$\widetilde{Z}^{ss}_{R'} = \{x \in \widetilde{X}^{ss,R}\,\,  | \,\, R' x =x \} = \widetilde{Z}_{R'} \cap X^{ss,R}$$
with $\widetilde{Z}_{R'}$ the union of those connected components of the fixed point set $\widetilde{X}^{R'}$ over which $R'$ acts trivially on the fibres of $\widetilde{\mathcal{L}}$. Moreover  $N_{R'}$ and its quotient group $N_{R'}/R'$ are reductive, and $\widetilde{Z}^{ss}_{R'} $ coincides with both the semistable locus and the stable locus for the induced action of $N_{R'}/R'$ on $\widetilde{Z}_{R'}$.

In this situation $R\widetilde{Z}^{ss}_{R'}  \setminus E \cong R \times^{N_{R'}} (\widetilde{Z}^{ss}_{R'}  \setminus E)$ is a nonempty $R$-invariant (and hence $H$-invariant) open subscheme of $\widetilde{X} \setminus E$ with a geometric quotient
$$(R\widetilde{Z}^{ss}_{R'}  \setminus E)/R \cong (\widetilde{Z}^{ss}_{R'}  \setminus E)/N_{R'}.$$
Since $\psi$ restricts to an $H$-equivariant isomorphism from $\widetilde{X} \setminus E$ to an open subscheme of $X$, it follows that if we set
$$\mathcal{S}_0(X,H,\lambda,\mathcal{L}) = \psi(R \widetilde{Z}^{ss}_{R'}  \setminus E)$$
then $\mathcal{S}_0(X,H,\lambda,\mathcal{L})$ is a nonempty  $R$-invariant (and hence $H$-invariant) open subscheme of $X$ with a geometric quotient
$$\mathcal{S}_0(X,H,\lambda,\mathcal{L})/H \cong (\widetilde{Z}^{ss}_{R'}  \setminus E)/N_{R'}.$$

\medskip

\noindent {\bf Case 2(c)ii).} By Remark \ref{rem:specialcases}, if we are not in the situation of 2(c)i), then the semistable locus $\widetilde{X}^{ss,R}$ is empty, and as in Case 2(b) there is a nonempty open stratum
$$\widetilde{S}_\beta \cong R \times^{P_\beta} \widetilde{Y}^{ss}_\beta$$
with $\beta \neq 0$
in the HKKN stratification of $\widetilde{X}$. Then $\widetilde{S}_\beta \setminus E \cong R \times^{P_\beta} (\widetilde{Y}^{ss}_\beta \setminus E)$ is nonempty and open in $\widetilde{X}$. 

If $\overline{\widetilde{Y}^{ss}_\beta \setminus E} \neq \widetilde{X}$ then as in Case 2(b)i)
$\mathcal{S}_0(X,H,\lambda,\mathcal{L}) = \psi(R \, \mathcal{S}_0(\overline{\widetilde{Y}^{ss}_\beta \setminus E},P_\beta,\lambda_\beta,\widetilde{\mathcal{L}}|_{\overline{\widetilde{Y}^{ss}_\beta \setminus E}}) \setminus E)$ is a nonempty $R$-invariant (and hence $H$-invariant) open subscheme of $X$ with a geometric quotient
$$\mathcal{S}_0(X,H,\lambda,\mathcal{L})/H \cong (\mathcal{S}_0(\overline{\widetilde{Y}^{ss}_\beta \setminus E},P_\beta,\lambda_\beta,\widetilde{\mathcal{L}}|_{\overline{\widetilde{Y}^{ss}_\beta \setminus E}}) \setminus E)/P_\beta.
$$
Finally if $\overline{\widetilde{Y}^{ss}_\beta \setminus E} = \widetilde{X}$ then as in Case 2(b)ii) we can either reduce to Case 1 or use induction on the dimension of $H$.

\medskip

This completes the recursive definition of $\mathcal{S}_0(X,H,\lambda,\mathcal{L})$; then we define the stratification  $\{ \mathcal{S}_\gamma | \gamma \in \Gamma \}$ and its indexing set $\Gamma$ as at (\ref{(1)}) and (\ref{(2)}), and Theorem \ref{mainthmstrat} follows from the construction. Of course, if the dimension of $X$ or $H$ is $0$,  we set $\mathcal{S}_0(X,H,\lambda,\mathcal{L}) = X$. 

\medskip

\noindent  {\bf Example: Ordered points on the projective line}.
This is a familiar example  (\textit{cf.}\ \cite{K,New}). Let $H=\SL(2)$ so that $U$ is trivial and $R=H$. Fix $n\ge 1$ and consider the diagonal action of $H$ on $(\PP^1)^n$. The linearisation is the $n$th tensor power of the standard representation on $\CC^2$ via the Segre embedding,  
and the indexing set of the HKKN stratification is $$\mathcal{B}=\{(2r-n) \mid n\ge r > n/2\}\cup \{0\}.$$
If $\b=2r-n$ with $r>n/2$ then a sequence lies in $Z_\b = Z_\b^{ss}$ if and only if it contains $r$ entries equal to $[1:0]$ and $n-r$ entries equal to  $[0:1]$, while $Y_\b= Y_\b^{ss}$ consists of sequences with precisely $r$ entries equal to  $[1:0]$, and finally the HKKN stratum $S_\b$ consists of sequences such that exactly $r$ entries coincide. Thus $Z_\b,Y_\b$ and $S_\b$ all have ${n \choose r}$ connected components. The semistable stratum corresponds to $\b=0$ and consists of sequences in which no point of $\PP^1$ occurs strictly more than $n/2$ times. 

In order to describe the refined stratification, note first  that the $U_\beta$-stabilisers in $Z_\beta$ when $\b=2r-n$ are trivial for $n/2<r<n$. Therefore in the refined stratification, $S_\b$ decomposes as the disjoint union of $S_{2r-n}^{r,n-r}=\SL(2)Z_\b$, consisting of sequences supported at two points of multiplicity $r$ and $n-r$ respectively, and its  complement $S_{2r-n}^{r,<n-r}$ in $S_\b$; when  $r=n-1$ this complement is empty. When $r=n$ the $U_n$-stabilisers are not trivial, but they have the same dimension for all points in $Z_{n}$ and the stratum $S_{n}$ remains intact.

The semistable stratum $S_0$ behaves differently for even and odd $n$. For odd $n$, the $\SL(2)$-semistable locus $S_0$ coincides with the stable locus and no blow-up is needed. However, for even $n$  we need to split $S_0$ into smaller strata to get geometric quotients: the stable locus $S_0^{<n}$ where strictly fewer than $n/2$ points coincide, the (connected components of the) centre $S_0^{n/2,n/2}$ of the blow-up where half the entries of the sequence coincide at each of two points in $\PP^1$, and the (connected components of the) strata $S_0^{n/2,<n/2}$ where half the entries coincide but the other half do not.

Thus modulo taking connected components the refined stratification has the following form:
\[\text{ n odd:   }(\PP^1)^n=S_0 \sqcup S_{n-2} \sqcup S_{n} \sqcup \bigsqcup_{n/2<r<n-1} \underbrace{S_{2r-n}^{r,n-r} \sqcup S_{2r-n}^{r,<n-r}}_{S_{2r-n}}; \]
\[\text{ n even:   }(\PP^1)^n=\underbrace{S_0^{<n} \sqcup S_0^{n/2,n/2} \sqcup S_0^{n/2,<n/2}}_{S_0} \sqcup S_{n-2} \sqcup S_{n} \sqcup \bigsqcup_{n/2<r<n-1} \underbrace{S_{2r-n}^{r,n-r} \sqcup S_{2r-n}^{r,<n-r}}_{S_{2r-n}}. \]

\begin{rmk}
If we only require good quotients rather than geometric quotients of the strata $S_\gamma$, then we can simplify the algorithm for constructing the stratification $\{ \mathcal{S}_\gamma | \gamma \in \Gamma \}$  somewhat, by combining the Cases 2(a) and (c) in the definition of $\mathcal{S}_0(X,H,\lambda,\mathcal{L})$ into the following single case.

\medskip

\noindent {\bf Case 2(a').} Suppose that the semistable locus $X^{ss,R}$ for the linear action of $R$ on $X$ is nonempty. Then we let 
$$\mathcal{S}_0(X,H,\lambda,\mathcal{L}) = X^{ss,R}$$
and by classical GIT this has a good quotient $X/\!/R$.
\end{rmk}

We note that this stratification is not intrinsic to the stack $\mathfrak{X}:=[X/H]$, as it depends on the presentation of $\mathfrak{X}$ as a quotient stack $[X/G]$ and a choice of ample linearisation $\cL$ on $X$, as well as an invariant inner product on $\Lie R$. Since for reductive groups $H=R$, this stratification refines the HKKN stratification, which itself depends on a choice of such a linearisation and an inner product, this dependence is to be expected. 

\begin{rmk}\label{strat proj}
Note that when $H=R$ is reductive then condition \eqref{dagger} in Theorem \ref{mainthmstrat} holds for $\beta = 0$ if and only if  $X^{ss} = X^s$. Moreover $(\dagger)$ holds for $\beta \neq 0$ if and only if semistability coincides with stability for the induced linear action of  $\text{Stab}(\beta)/\lambda_\beta(\GG_m)$ on $Z_\beta$ and also condition \eqref{star} holds for the action of the graded unipotent radical $\hU_\beta = U_\beta \rtimes \lambda_\beta(\GG_m)$ of the parabolic subgroup $P_\beta \leq R$.
\end{rmk}

There is a notion of a coarse moduli space of an Artin stack $\mathfrak{X}$ (over $\kk$), which is a morphism $\pi : \mathfrak{X} \ra Y$ to an algebraic space $Y$ that is a categorical quotient (i.e. $\pi$ is initial among all morphisms from $\mathfrak{X}$ to an algebraic space) and such that $\pi$ induces a bijection on point sets $[\mathfrak{X}(\kk)] \cong Y(\kk)$. 

\begin{rmk}\label{geom implies cms}
In general, if $G$ is an algebraic group acting on a scheme $X$ with a categorical quotient $q : X \ra Y$, then $\pi : [X/G] \ra Y$ is not necessarily a coarse moduli space, as $\pi$ only induces a bijection on points sets if $Y$ is an orbit space. However, if $q$ is a geometric quotient, then $\pi$ is a coarse moduli space for $\mathfrak{X}:=[X/G]$. Indeed, since geometric quotients are orbit spaces, $\pi$ clearly induces a bijection $[\mathfrak{X}(\kk)]  = X(\kk)/G(\kk) \cong Y(\kk)$. To prove it is also a categorical quotient, one can apply a result of Rydh \cite[Theorem 3.16]{Rydh}, which says that if $q$ is universally open and strongly geometric (in the sense of \cite[Definition 2.2]{Rydh}), then $\pi$ is a categorical quotient. Since $q$ is a geometric quotient it is surjective, universally open and the fibres of $q$ are single orbits; furthermore, the natural morphism $j : G \times X \ra X \times_Y X$ is surjective. In particular, $j$, $q$ and $q^{\mathrm{cons}}$ (the induced map for the constructible topology) are all universally submersive; thus $q$ is a strongly geometric quotient by \cite[Lemma 2.3]{Rydh}.
\end{rmk}

Alper \cite{alper} introduced a notion of a good moduli space for stacks, which is inspired by the notion of a good quotient in reductive GIT. However, this notion is suited to reductive groups and  is less favourable for our set up where we have stacks with potentially non-reductive stabiliser groups. Indeed  $BG = [\spec \kk /G] \ra \spec \kk$ is a good moduli space if and only if $G/\kk$ is linearly reductive, as the stabiliser groups of the closed points of a stack admitting a good moduli space are all linearly reductive by \cite[Proposition 12.14]{alper}.

By Remark \ref{geom implies cms}, we obtain the following corollary of Theorem \ref{mainthmstrat}.

\begin{corollary}\label{cor cms}
Let $\mathfrak{X}= [X/H]$ be a quotient stack, where $X$ and $H$ are as in the statement of Theorem \ref{mainthmstrat}, and so we have a stratification $\{ \mathcal{S}_\gamma | \gamma \in \Gamma \}$ of $X$ into $H$-invariant quasi-projective subschemes which each admit geometric quotients $S_\gamma / H$. Then there is an induced stratification $\{ \mathfrak{X}_\gamma := [\mathcal{S}_\gamma /H] : \gamma \in \Gamma \}$ of $\mathfrak{X}$ such that each stratum has a coarse moduli space $\mathfrak{X}_\gamma \ra \mathcal{S}_\gamma /H$.
\end{corollary}

We note that, in general, it is not the case that $ [\mathcal{S}_\gamma /H]$ is isomorphic to $\mathcal{S}_\gamma /H$; this only happens if additionally the action of $H$ on $\mathcal{S}_\gamma$ is free.

\section{Stacks filtered by quotient stacks}\label{sec filt stack}

In $\S$\ref{sec strat quot stacks}, we constructed a stratification $\{ \mathcal{S}_\gamma | \gamma \in \Gamma \}$ of a projective scheme $X$ acted on linearly by a linear algebraic group $H=U\rtimes R$ with internally graded unipotent radical $U$,  such that each stratum is an $H$-invariant quasi-projective subscheme of $X$ with a geometric quotient by the action of $H$. The quotient stack $[X/H]$ then has an induced stratification $\{\Sigma_\gamma | \gamma \in \Gamma \}$ such that each stratum $\Sigma_\gamma$ has the form 
$$ \Sigma_\gamma := [\mathcal{S}_\gamma /H] \cong [W_\gamma/H_\gamma]$$
where $W_\gamma$ is a quasi-projective subscheme of $X$ acted on by a linear algebraic subgroup $H_\gamma$ of $H$ with internally graded unipotent radical, and this action has a geometric quotient $W_\gamma/H_\gamma$. Moreover under appropriate hypotheses of \lq semistability coinciding with stability' for the actions of the subgroups $H_\gamma$, the geometric quotients $W_\gamma/H_\gamma$ are themselves projective.

We can apply this to Artin stacks which are filtered by quotient stacks admitting stratifications that are compatible in the following sense. 

\begin{definition}\label{def filt stack}
Let $\cM$ be an Artin stack.
\begin{enumerate}
\item We say $\cM$ is {\em filtered by quotient stacks} if $\cM$ can be expressed as an increasing union 
$$\mathcal{M} = \bigcup_{n \geqslant 0} \cM_n $$
of open substacks with presentations of the form
$$ \cM_n \cong [V_n/H_n]$$
where $[V_n/H_n]$ is the quotient stack associated to a linear action on a quasi-projective scheme $V_n$ by a group $H_n = U_n \rtimes R_n$ with internally graded unipotent radical. 
\item We say $\cM$ is {\em filtered by quotient stacks admitting compatible stratifications} if $\cM$ is filtered as above and moreover, for each $n \in \mathbb{N}$, there is an ample $H_n$-linearisation $\cL_n$ on projective completion $\overline{V_n}$ of $V_n$ and a choice of invariant inner product on the Lie algebra of the Levi factor $R_n \leq H_n$ such that for $n' > n$, the pullback along $\cM_n \hookrightarrow \cM_{n'}$ of the stratification of $\cM_{n'}$ (obtained by restricting the stratification on $[\overline{V_{n'}}/H_{n'}]$ defined by $\cL_{n'}$ and this norm to $\cM_{n'}$) coincides with the corresponding stratification on $\cM_n$.
\item We additionally say that $\cM$ is {\em filtered by quotient stacks admitting compatible stratifications which are asymptotically stable}, if for each $n$ and each stratum $\Sigma_{\gamma,n}$ of $\cM_n \cong [V_n/H_n]$ there exists $N \geq n$ such that if $l>m\geqslant N$ and $\Sigma_{\gamma,l} \hookrightarrow \cM_l$ and $\Sigma_{\gamma,m} \hookrightarrow \cM_m$ are the strata whose restrictions to $\cM_n$ are $\Sigma_{\gamma,n}$, then $\Sigma_{\gamma,l}$ is the image of $\Sigma_{\gamma,m}$ under the inclusion $\cM_m \hookrightarrow \cM_l$.
\end{enumerate}
\end{definition}

Then we obtain the following corollary to Theorem \ref{mainthmstrat}.

\begin{corollary}\label{strat filt stacks}
Let $\cM$ be an Artin stack that is filtered by quotient admitting compatible stratifications which are asymptotically stable; then there is a stratification $\{\Sigma_\gamma | \gamma \in \Gamma \}$ of the stack $\mathcal{M}$ such that each stratum $\Sigma_\gamma$ is isomorphic to a quotient stack 
$ [W_\gamma/H_\gamma]$,
where $W_\gamma$ is quasi-projective and acted on by a linear algebraic group $H_\gamma = U_\gamma \rtimes R_\gamma$  with internally graded unipotent radical $U_\gamma$. Moreover, there is a geometric quotient $W_\gamma/H_\gamma$ which is a coarse moduli space for $\Sigma_\gamma$. In particular, $\cM$ has a stratification by quotient stacks which all have coarse moduli spaces.
\end{corollary}

\begin{rmk}
In fact, the same result holds if we consider stacks which admit an increasing open cover indexed by a partially ordered filtered set that contains a cofinal copy of $\mathbb{N}$ with analogous properties.
\end{rmk}

We have the following criterion for these moduli spaces to be projective (\textit{cf.}\  Remark \ref{strat proj}).

\begin{rmk}
In the situation of Corollary \ref{strat filt stacks}, for a stratum $\Sigma_\gamma \cong  [W_\gamma/H_\gamma]$ where the unipotent radical $U_\gamma$ is graded by $\lambda_\gamma:\GG_m \to R_\gamma$, the coarse moduli space $W_\gamma/H_\gamma$ will be projective if \eqref{star} holds for the action of $U_\gamma \rtimes \lambda_\gamma(\GG_m)$ on a suitable projective completion of $W_\gamma$ with respect to an induced linearisation, and semistabiity coincides with stability for the induced action of $R_\gamma/\lambda_\gamma(\GG_m)$ on the analogue of $Z_{\min}$ in this projective completion.
\end{rmk}

\section{Applications}\label{sec appl}

We can apply the previous section to  the construction of moduli spaces of \lq unstable' objects, including moduli of sheaves over a projective scheme \cite{behjk16}, projective curves \cite{jj}) and hypersurfaces in a complete toric variety.

\subsection{Moduli of hypersurfaces in a complete toric variety}

Let $Y$ be a complete simplicial toric variety. It was observed in \cite[$\S$4]{BDHK2} (using the description of $\Aut(Y)$ given in \cite{cox}) that the automorphism group $H=\Aut(Y)$ of $Y$ is a linear algebraic group with internally graded unipotent radical.

By considering suitable Hilbert schemes, each component of the moduli stack of hypersurfaces in $Y$ can be expressed itself as a quotient stack $[V/H]$ where $V$ is a projective scheme with a linear action of $H:=\Aut(Y)$ with an ample linearisation. Thus we can apply Theorem \ref{mainthmstrat} directly to the moduli stack  $\mathcal{M}$ to obtain a stratification  $\{\Sigma_\gamma | \gamma \in \Gamma \}$ of $\mathcal{M}$ such that each stratum $\Sigma_\gamma$ is isomorphic to a quotient stack 
$ [W_\gamma/H_\gamma]$,
where $W_\gamma$ is a quasi-projective scheme acted on by a linear algebraic group $H_\gamma$  with internally graded unipotent radical, and there is a geometric quotient $W_\gamma/H_\gamma$ which is a coarse moduli space for $\Sigma_\gamma$.

A simple example is the automorphism group of the weighted projective plane $Y=\PP(1,1,2) = (\kk^3 \setminus \{ 0 \})/ \GG_m$, for $\GG_m$ acting linearly on $\kk^3$ with weights $1,1,2$. This is given by 
$$\Aut(\PP(1,1,2)) \cong U \rtimes R$$
where $R \cong \GL(2)$ is reductive and 
$U \cong (\GG_a)^3$  is unipotent, with elements  $(\lambda,\mu,\nu) \in \kk^3$ acting on $\PP(1,1,2)$ via 
$$[x,y,z]  \mapsto [x,y,z+\lambda x^2 + \mu xy + \nu y^2].$$
The central one-parameter subgroup $\GG_m$ of $R \cong \GL(2)$ acts on the Lie algebra of 
$U$ with all positive weights, providing an internal  grading of the unipotent radical $U$. 

\subsection{Moduli of sheaves over a projective scheme} 

Let $(Y,\cO_Y(1))$ denote a polarised projective scheme over an algebraically closed field $\kk$ of characteristic $0$. For a fixed Hilbert polynomial $P$, we let 

$\cM:= \cM_{Y,P}$ denote the moduli stack of coherent sheaves over $Y$ with Hilbert polynomial $P$ with respect to $\cO_Y(1)$. We can filter this stack by quotient stacks in several ways: for example, one can use Castelnuevo--Mumford regularity or the Harder--Narasimhan type. In both filtrations, we will filter by quotient stacks which are quotients by general linear groups of open subschemes of Quot schemes, which appear in Simpson's GIT construction \cite{s94} of the moduli space of semistable sheaves on $Y$ with Hilbert polynomial $P$. In this section, we will describe the two ways to filter $\cM$ and explain how to use a filtration involving Harder--Narasimhan types to apply the results in $\S$\ref{sec filt stack}.

Let us recall Simpson's reductive GIT construction of moduli spaces of semistable sheaves. By the Le Potier estimates \cite[Theorem 1.1]{s94}, the family of semistable sheaves over $Y$ with Hilbert polynomial $P$ is bounded; hence, for $n \gg 0$, all such sheaves are $n$-regular in the sense of Castelnuevo-Mumford. For an $n$-regular sheaf $\cE$, the evaluation map 
\[ H^0(\cE(n)) \otimes \cO_Y(-n) \ra \cE\]
is surjective and $\dim H^0(\cE(n)) = P(\cE,n)$, as all the higher cohomology groups of $\cE(n)$ vanish; hence, any $n$-regular sheaf with Hilbert polynomial $P$ determines a closed point in the Quot scheme 
\[Q_n:=\quot(\cO_Y(-n)^{\oplus P(n)},P)\] 
parameterising quotient sheaves of $\cO_Y(-n)^{\oplus P(n)}$ with Hilbert polynomial $P$. In fact, it determines a point in the open subscheme $Q_{n-\reg}$ consisting of quotients $q : \cO_Y(-n)^{\oplus P(n)} \twoheadrightarrow \cF$ such that $\cF$ is $n$-regular and $H^0(q(n))$ is an isomorphism. If $Q_{n-\reg}^{ss}$ denotes the open subset, where additionally the quotient sheaf is semistable, then Simpson constructs the moduli space of semistable sheaves as a quotient of the natural $\SL_{P(n)}$-action on the closure of $Q_{n-\reg}^{ss}$. The action is linearised by choosing $m \gg n$ to construct an embedding of the Quot scheme into a Grassmannian and then composing with a Pl\"ucker embedding of the Grassmannian into projective space; we let $\mathcal{L}_{m,n}$ denote this linearisation. 

Since $n$-regularity is an open condition, we obtain an increasing open filtration of $\cM$ by the substacks $\cM_{n-\reg}$ of $n$-regular sheaves; furthermore, there are isomorphisms $\cM_{n-\reg} \cong [Q_{n-\reg}/\GL_{P(n)}]$ for all $n$. Hence $\cM$ is filtered by quotient stacks using $n$-regularity. 

Let us explain a second way to filter $\cM$ using Harder--Narasimhan (HN) types \cite{HN}. We recall that every coherent sheaf $\cE$ has a unique HN filtration $ 0 = \cE^{(0)} \subsetneq \cE^{(1)} \cdots \cE^{(r)}= \cE$ (where here we use a Rudakov type notion of semistability to define HN filtrations for any coherent sheaf; for example, see \cite[Definition 4.11]{H}). The HN type of $\cE$ is $\tau(\cE) :=(P(\cE_1), \dots , P(\cE_r))$, where $\cE_i := \cE^{(i)}/ \cE^{(i-1)}$. The set of HN types for sheaves with Hilbert polynomial $P$ is partially ordered, such that the semistable HN type given by the single polynomial $(P)$ is minimal. We let $\cM^\tau$ and $\cM^{< \tau}$ denote the substacks of sheaves of HN type $\tau$ and of HN type less than $\tau$ respectively. By work of Shatz \cite{shatz}, the HN type is upper semi-continuous, and so $\cM^{< \tau}$ is an open substack of $\cM$. Furthermore, as there are only finitely many HN types less than a given HN type $\tau$, these sheaves form a bounded family, and so $\cM^{<\tau}$ is a quotient stack of the form $[Q_{n-\reg}^{<\tau}/\GL_{P(n)}]$ for $n \gg 0$ (depending on $\tau$); here, $Q_{n-\reg}^{<\tau}$ is the open subscheme of $Q_{n-\reg}$ parametrising quotient sheaves with HN type less than $\tau$.

\begin{rmk}\label{compare HN and HKKN}
The relationship between the HKKN stratifications and the HN stratifications on $Q_{n-\reg}$ is studied in \cite{H,hok12}, where it is shown that for a HN type $\tau$ and for $m \gg n \gg 1$ (depending on $\tau$), there is an HKKN index $\beta_{n,m}(\tau)$ for the $\SL_{P(n)}$-action on $Q_n$ with respect to $\cL_{n,m}$ (and the Killing form on $\SL_{P(n)}$) such that the following statements hold.
\begin{enumerate}
\item All sheaves with HN type $\tau$ are $n$-regular, and thus $\cM_\tau \cong [Q_{n-\reg}^\tau/\GL_{P(n)}]$, where $Q_{n-\reg}^\tau$ is the HN stratum in $Q_{n-\reg}$ indexed by $\tau$.
\item The HN stratum $Q_{n-\reg}^\tau$ is contained as a closed subscheme in the HKKN stratum $S_{\beta_{n,m}(\tau)}$.
\end{enumerate}
For the precise definition of $\beta_{n,m}(\tau)$, see \cite[Definition 4.18]{H}. For fixed $n$ and $m$, the assignment $\tau \mapsto \beta_{n,m}(\tau)$ is not injective; however, for two HN types $\tau \neq \tau'$ and $m \gg n \gg 0$, we have  $\beta_{n,m}(\tau)\neq \beta_{n,m}(\tau')$. In fact, for a finite number of HN types $\tau_1, \dots, \tau_l$, we can pick $m \gg n \gg 1$ so that all the indices $\beta_{n,m}(\tau_i)$ are distinct and both the  properties (1) and (2) above hold for all of these HN types. 

Furthermore, in \cite[Theorem 1.1]{H} an asymptotic HKKN stratification on $\cM$ is constructed using the HKKN stratifications on each $\cM_{n-\reg}$ and it is shown that this asymptotic HKKN stratification coincides with the HN stratification on $\cM$. In fact, the HKKN stratifications on $\cM_{n-\reg}$ are not compatible stratifications which are asymptotically stable in the sense of Definition \ref{def filt stack}, and so this asymptotic construction of a HKKN stratification is slightly more involved (for details, see \cite[$\S$4]{H}).
\end{rmk}

We claim that if we consider the filtration given by the open substacks $\cM^{<\tau}$ and choose an appropriate presentation of these stacks, then the HKKN stratifications on each $\cM^{<\tau}$ coincide with the HN stratifications, and thus in particular, these stratifications are compatible and asymptotically stable in the sense of Definition \ref{def filt stack}. For each HN type $\tau$, as there are only finitely many HN types less than $ \tau$, we can pick $m \gg n \gg 1$ as in Remark \ref{compare HN and HKKN}, so that both statements in this remark hold for all these HN types and they all have different HKKN indices. Indeed, if we compare the stratifications on $Q_{n-\reg}^{<\tau}$, then we have $Q_{n-\reg}^{\tau'} \hookrightarrow S_{\beta_{n,m}(\tau')}$ for all $\tau' < \tau$, and as $Q_{n-\reg}^{<\tau}$ is the disjoint union of the HN strata $Q_{n-\reg}^{\tau'}$, we see that in fact the HN and the HKKN stratifications must agree.

Finally let us consider the refinement of the HKKN stratification on each $\cM^{<\tau}$ provided by Theorem \ref{mainthmstrat}. Since the refinement will come from considering various stabiliser groups and stable loci, and these can be studied intrinsically  for the sheaves themselves, we see that the stratifications on each $\cM^{<\tau} \cong [Q_{n-\reg}^{<\tau}/\GL_{P(n)}]$ are compatible and asymptotically stable in the sense of Definition \ref{def filt stack}. Therefore, there is an induced stratification  $\{\Sigma_\gamma | \gamma \in \Gamma \}$ of $\cM$ by quotient stacks that refines the HN stratification and such that each quotient stack has a presentation $\Sigma_\gamma \cong  [W_\gamma/H_\gamma]$, where $W_\gamma$ is a quasi-projective scheme acted on by a linear algebraic group $H_\gamma$  with internally graded unipotent radical, and there is a geometric quotient $W_\gamma/H_\gamma$ which is a coarse moduli space for $\Sigma_\gamma$. In fact, we can say what the open stratum is: provided that there exists a stable sheaf with Hilbert polynomial $P$, the open stratum is the moduli stack $\cM^{s}$ of stable sheaves. The precise construction and description of these coarse moduli spaces is described in \cite{behjk16}.

\subsection{Moduli of unstable curves}

Another classical family of moduli spaces constructed as reductive GIT quotients is given by the moduli spaces of stable curves of fixed genus $g \geq 2$. Here  we can use  special linear group actions on Hilbert schemes of curves embedded in projective spaces, and exploit the Rosenlicht--Serre description of singular curves in terms of their normalisations together with additional data describing the singularities, to 
 look for an associated stratification $\{\Sigma_\gamma | \gamma \in \Gamma \}$ of the moduli stack of projective curves, allowing us to  construct geometric quotients of the unstable strata \cite{jj}. In this case of projective curves the linear actions can be set up so that (almost always) the condition that semistability coincides with stability  is satisfied. Thus we expect  the coarse moduli spaces of \lq unstable curves' given by the  geometric quotients $W_\gamma/H_\gamma$ of the strata $\Sigma_\gamma \cong [W_\gamma/H_\gamma]$ to be projective.

\end{document}